\documentclass[article]{siamltex}
\usepackage[english]{babel}
\usepackage[labelformat=simple,font={it}]{subcaption}   
\usepackage[font={small,it}]{caption}
\usepackage{amsmath,amssymb,amsfonts,graphicx,outlines}
\usepackage{myshortcuts}
\usepackage{xcolor,bm}
\pagecolor{white}
\usepackage{fancyhdr}
\usepackage{MnSymbol}
\usepackage{csvsimple,mleftright}
\usepackage{tikz}
\usepackage{pgfplots}
\usepackage[normalem]{ulem}
\newcommand{\stkout}[1]{\ifmmode\text{\sout{\ensuremath{#1}}}\else\sout{#1}\fi}
\usetikzlibrary{decorations.pathreplacing}

\DeclarePairedDelimiter\abs{\lvert}{\rvert}%

\usepackage{appendix}
\usepackage{comment}
\usepackage{relsize}
\newcommand{\norm}[1]{\left\lVert#1\right\rVert}
\usepackage[linesnumbered,lined,boxed,ruled,vlined]{algorithm2e}
\let\oldnl\nl
\newcommand{\nonl}{\renewcommand{\nl}{\let\nl\oldnl}}

\newcommand*\samethanks[1][\value{footnote}]{\footnotemark[#1]}
\newtheorem{remark}[theorem]{Remark}

\newenvironment{customlegend}[1][]{%
    \begingroup
    \csname pgfplots@init@cleared@structures\endcsname
    \pgfplotsset{#1}%
}{%
    \csname pgfplots@createlegend\endcsname
    \endgroup
}%
\def\addlegendimage{\csname pgfplots@addlegendimage\endcsname}

\title{Preconditioned Chebyshev BiCG for parameterized linear systems}

\author{Siobhán Correnty\thanks{Department of Mathematics, Royal Institute of Technology (KTH), SeRC Swedish
e-Science Research Center, Lindstedtsvägen 25, Stockholm, Sweden
\{\texttt{correnty,eliasj}\}\texttt{@kth.se}} \and
Elias Jarlebring\samethanks \and Daniel B. Szyld\thanks{Department of Mathematics, Temple University, 1805 N. Broad Street, Philadelphia, PA 19122-6094, \texttt{szyld@math.temple.edu}} }

\begin{document}

\maketitle

\begin{abstract}
We consider the problem of approximating the solution to $A(\mu) x(\mu) = b$ for many different values of the parameter $\mu$. Here $A(\mu)$ is large, sparse, and nonsingular with a nonlinear dependence on $\mu$. Our method is based on a companion linearization derived from an accurate Chebyshev interpolation of $A(\mu)$ on the interval $[-a,a]$, $a \in \mathbb{R}_+$, inspired by Effenberger and Kressner [BIT, 52 (2012), pp. 933–951]. The solution to the linearization is approximated in a preconditioned BiCG setting for shifted systems, as proposed in Ahmad et al. [SIAM J. Matrix Anal. Appl., 38 (2017), pp. 401–424], where the Krylov basis matrix is formed once. This process leads to a short-term recurrence method, where one execution of the algorithm produces the approximation to $x(\mu)$ for many different values of the parameter $\mu \in [-a,a]$ simultaneously. In particular, this work proposes one algorithm which applies a shift-and-invert preconditioner exactly as well as an algorithm which applies the preconditioner inexactly based on the work by Vogel [Appl. Math. Comput., 188 (2007), pp. 226–233]. The competitiveness of the algorithms is illustrated with large-scale problems arising from a finite element discretization of a Helmholtz equation with parameterized material coefficient. The software used in the simulations is publicly available online, and thus all our experiments are reproducible. 
\end{abstract}  




\section{Introduction}
In this work, we propose a new approach for computing an efficient approximation to the solution of the parameterized linear system given by
\begin{align} \label{our-problem}
	A(\mu) x(\mu) = b
\end{align}
for many values of the parameter $\mu$ simultaneously. Here $A(\mu) \in \mathbb{R}^{n \times n}$ is assumed nonsingular, nonlinear in $\mu \in [-a,a]$, $a \in \mathbb{R}_+$, and $b \in \mathbb{R}^n$. Specifically, our method finds accurate approximations of $x(\mu)$ for $\mu$ in a specified region of the interval $[-a,a]$. Parameterized linear systems have been studied previously, for example, in \cite{FrommerMaass99:20,KilmerOleary01}, where these systems arise in the context of Tikhonov regularization for ill-posed problems, as well as in \cite{KressnerToblerLR}, where the solution was approximated by a tensor of low rank, and in \cite{GuSimoncini}, where the right-hand side vector also depended on the parameter. 

We assume further that $A(\mu)$ in \eqref{our-problem} is large and sparse, and can be expressed as the sum of products of matrices and functions, i.e., 
\begin{align} \label{summation}
	A(\mu) = C_1 f_1 (\mu) + \cdots + C_{n_f} f_{n_f}(\mu), 
\end{align}
where $n_f \ll n$. Our method requires an approximation of $A(\mu)$ via a Chebyshev interpolation. In this way, we compute $P(\mu) \approx A(\mu)$, where
\begin{align} \label{chebyshev}
	P(\mu) =  P_0 \tau_0 (\mu) + \ldots + P_d \tau_d (\mu) 
\end{align}
with $P_{\ell} \in \mathbb{R}^{n \times n}$ and $\tau_{\ell} (\mu)$ the recursively defined Chebyshev polynomials on the interval $[-a,a]$. The matrix $P(\mu)$ is assumed nonsingular throughout this work, and we assume that $d$, the truncation parameter in the Chebyshev approximation \eqref{chebyshev}, is large enough such that the error introduced by the Chebyshev interpolation is small. 

We present a preconditioned short-term recurrence Krylov subspace method to approximate the solution to 
\begin{align} \label{approx}
	P(\mu) \tilde{x}(\mu) = b,
\end{align}
where $\tilde{x}(\mu) \approx x(\mu)$. In practice, our method approximates the solution to a companion linearization of the form
\begin{align} \label{companion}
	(K - \mu M ) u(\mu) = \tilde{b}
\end{align}
with coefficient matrices $K$, $M \in \mathbb{R}^{dn \times dn}$ and constant vector $\tilde{b} \in \mathbb{R}^{dn}$. The solution to \eqref{approx} and the companion linearization are equivalent in a certain sense, shown in Section~\ref{sec:linearization}. Here the bases for two Krylov subspaces are generated via a Lanczos biorthogonalization procedure as in the biconjugate gradient method (BiCG) \cite{BiCG76,Lanczos52}. Specifically, the method incorporates shift-and-invert preconditioners of the form $(K-\sigma M)^{-1}$ and $(K-\sigma M)^{-T}$, for $\sigma \in (-a,a)$, to accelerate convergence for solutions corresponding to values of $\mu$ close to the chosen target $\sigma$. Additionally, the use of such preconditioners leads to a shifted linear system, and shift- and scaling-invariance properties of Krylov subspaces are exploited. In this way, we reuse one Krylov subspace basis matrix to determine approximations to \eqref{our-problem} for many different values of $\mu$. 

We propose two variants of our method. The first variant considers an exact application of the preconditioner in a BiCG setting for shifted systems, inspired by \cite{PreCondAhmed}. The second variant incorporates an approximate application of the preconditioners in an inexact BiCG setting for shifted systems, based on the prior works \cite{Saad1993,Vogel2007}. We prove a bound on the residual of the second variant, and the bound is efficient to compute.

The first variant is appropriate only when an LU decomposition of a matrix of dimension $n \times n$ is feasible, whereas the inexact variant has the potential to solve a wider variety of large-scale problems. While the second variant is based on the first variant, the first variant may be useful in itself in cases where the exact LU decomposition is computable, as discussed in Section~\ref{sec:Flex}. Numerical simulations from time-delay systems and a parameterized Helmholtz equation show the performance of our proposed algorithms. Note, BiCG with right preconditioning is used throughout this work. A left preconditioned setting would have been possible with the first variant but not the second, as inexact preconditioning requires right preconditioning.

This paper is organized as follows. In Section~\ref{sec:linearization} we present the Chebyshev linearization, and, in Section~\ref{sec:Shifted-Systems-BiCG}, we describe how an equivalent shifted linear system is obtained. This section also provides preliminaries for the method BiCG for shifted systems. Section~\ref{sec:Preconditioning} shows how the preconditioners are implemented in an efficient manner, and, in Section~\ref{sec:Deriv-MCBiCG}, we derive the preconditioned Chebyshev BiCG method for parameterized linear systems. Section~\ref{sec:Helmholtz} provides a numerical example from a discretized Helmholtz equation, and Section~\ref{sec:delay-problem} utilizes our method for computing the transfer function from a time-delay system. In Section~\ref{sec:Flex} we derive the inexact variant of the method. Furthermore, we prove a bound on the residual produced by iterates of the inexact method. This section also highlights the performance of our approach for solving large-scale parameterized systems effectively. Conclusions are given in Section~\ref{sec:conclusion}.

\section{Linearization} \label{sec:linearization}
We consider a technique called companion linearization, used in prior works on polynomial eigenvalue problems \cite{Mackey:2006:VECT}, as well as in works for parameterized linear systems; see, e.g., \cite{GuSimoncini}. Our proposed linearization is of the form \eqref{companion}. The solution to~\eqref{companion} and $\tilde{x}(\mu)$ in \eqref{approx} with $P(\mu)$ \eqref{chebyshev} are equivalent in a certain sense, described as follows. 

On the interval $[-a,a]$, the Chebyshev polynomials are defined by the well-known three-term recurrence 
\begin{subequations} \label{cheby-rec}
\begin{eqnarray} 
	\tau_0 (\mu) &\coloneqq& 1, \\
	\tau_1 (\mu) &\coloneqq& \frac{1}{a} \mu, \\
	\tau_{\ell+1}(\mu) &\coloneqq& \frac{2}{a} \mu \tau_{\ell}(\mu) - \tau_{\ell-1}(\mu),
\end{eqnarray}
\end{subequations}
and the interpolation condition 
\begin{align} \label{interpolation-condition}
	P(\mu_{\ell}^*) = A(\mu_{\ell}^*), \quad \ell=1,\ldots,d,
\end{align}
holds, where $\mu_{\ell}^*$ are the $d$ roots of the degree $d$ polynomial $\tau_{d}$. The companion linearization, adapted from the work \cite{ChebyEffenKress} and utilized here, is given by
\begin{align}\label{linear} 
\left(\begin{bmatrix}
0 \hspace{-.15cm}	& I  \hspace{-.15cm} & 				& 				& 				& \\
I  \hspace{-.15cm}	& 0 \hspace{-.15cm}	& I \hspace{-.15cm}  &	 			&	 			& \\
				& I \hspace{-.15cm}	& 0 \hspace{-.15cm}  & I \hspace{-.15cm}	& 				& \\
				& 				& 				&\ddots \hspace{-.15cm}&  			& \\
				& 				& 				& I \hspace{-.15cm}	& 0 \hspace{-.18cm}& I \\
P_0 \hspace{-.15cm}& P_1 \hspace{-.15cm}& \cdots	\hspace{-.15cm}& P_{d-3} \hspace{-.15cm}& (-P_d \hspace{-.04cm}+ \hspace{-.04cm}P_{d-2}) \hspace{-.15cm}& P_{d-1} 
\end{bmatrix}
\hspace{-.03cm}-\hspace{-.03cm}  \frac{\mu}{a}
\begin{bmatrix}
 I \hspace{-.18cm}& & & & & \\
& 2 I\hspace{-.18cm}& & & & \\
& & 2 I \hspace{-.18cm}& & & \\
& & & \ddots \hspace{-.18cm}& & \\
& & & & 2 I \hspace{-.18cm}& \\
& & & & & \hspace{-.18cm}-2 P_d \\
\end{bmatrix} \right)
\begin{bmatrix}
\hspace{-.08cm}u_0(\mu)\hspace{-.08cm} \\
\hspace{-.08cm}u_1(\mu)\hspace{-.08cm} \\
\hspace{-.08cm}u_2(\mu)\hspace{-.08cm} \\
\hspace{-.08cm}\vdots\hspace{-.08cm} \\
\hspace{-.08cm}u_{d-2}(\mu)\hspace{-.08cm} \\
\hspace{-.08cm}u_{d-1}(\mu)\hspace{-.08cm}
\end{bmatrix}
\hspace{-.03cm}  =\hspace{-.03cm}  
\begin{bmatrix}
0 \\
0 \\
0 \\
\vdots \\
0 \\
b 
\end{bmatrix}.
\end{align}
Specifically, $u_{\ell}(\mu) \coloneqq \tau_{\ell}(\mu) \tilde{x}(\mu) \in \RR^n$, for $\ell = 0,\ldots,d-1$, where $\tilde{x}(\mu)$ is the unique solution in \eqref{approx} and 
\begin{align} \label{b-tilde}
	\tilde{b} \coloneqq \begin{bmatrix} 0 & 0 & 0 & \cdots & 0 & b 
\end{bmatrix}^T \in \mathbb{R}^{dn}. 
\end{align}
In this way, \eqref{linear} is of the form described in \eqref{companion}, where we have made the substitution 
\begin{align} \label{recur-sub}
	P_d u_d (\mu) = P_d \left( \frac{2}{a} \mu u_{d-1}(\mu) - u_{d-2}(\mu) \right)
\end{align}
in the last block row, using the recurrence relation \eqref{cheby-rec}. Note that since $A(\mu)$ is as in~\eqref{summation}, the coefficient matrices $P_{\ell}$ used in the linearization can be computed efficiently using a discrete cosine transform \cite{DCT} of the scalar functions $f_i$, $i=1,\ldots,n_f$. Specifically, we compute $P_{\ell} = C_1 p_{\ell}^1 + \ldots + C_{n_f} p_{\ell}^{n_f}$, for $\ell = 0,\ldots,d$, where the $i$th function in \eqref{summation} is approximated by a degree $d$ polynomial, i.e., $f_i(\mu) \approx p_{0}^i \tau_0 (\mu) + \ldots + p_{d}^i \tau_{d} (\mu)$. The interpolations are performed using \texttt{Chebfun} \cite{ChebfunGuide} in Matlab, and the coefficients $p_{\ell}^i$ smaller in magnitude than $10^{-16}$ are set to zero.

Thus, we can consider highly accurate interpolations of $A(\mu)$ without substantial computation. The following theorem describes the equivalence of solutions of the approximation in~\eqref{approx} and the system \eqref{linear}, where the parameter $\mu$ appears only linearly. 

\begin{theorem} \label{theorem1}
Let $A(\mu)$ be as in \eqref{our-problem} with parameter $\mu \in \RR$ and $P(\mu)$ as in \eqref{chebyshev} such that \eqref{approx} has a unique solution $\tilde{x}(\mu)$. Then the linear system \eqref{linear} has a unique solution of the form $u(\mu) = \begin{bmatrix} u_0(\mu), \ldots, u_{d-1}(\mu) \end{bmatrix}^T$ with
\begin{align} \label{unique-sol}
	u_{\ell}(\mu) \coloneqq \tau_{\ell} (\mu) \tilde{x}(\mu), \quad \ell = 0,\ldots,d-1. 
\end{align}
\end{theorem}
{\em Proof.}
Consider the first $d-1$ block rows of \eqref{linear}. An induction using the three-term recurrence of the Chebyshev polynomials on the interval $[-a,a]$ as in \eqref{cheby-rec} implies $u_{\ell}(\mu) = \tau_{\ell}(\mu) \tilde{z}$, $\ell=0,\ldots,d-1$, for $\tilde{z} \in \RR^{n}$. Inserting this parameterized solution into the last block row of~\eqref{linear} gives
\begin{align*}
	 \Big( P_0 \tau_0(\mu) + \cdots + P_{d-3} \tau_{d-3}(\mu) + (-P_d + P_{d-2}) \tau_{d-2}(\mu)  + \frac{2}{a}\mu P_d \tau_{d-1}(\mu) \Big) \tilde{z}
	&= P(\mu) \tilde{z} \\ &= b,
\end{align*}
due to \eqref{chebyshev}. Note the substitution $P_d \tau_d(\mu) \tilde{z} = P_d \big( (2\mu /a) \tau_{d-1} (\mu)- \tau_{d-2}(\mu) \big) \tilde{z}$, analogous to the relation \eqref{recur-sub}. Thus, $\tilde{z} = \tilde{x}(\mu)$, and the solution \eqref{unique-sol} to the linear system \eqref{linear} is unique since $\tilde{x} (\mu)$ is unique. $\endproof$

\begin{remark}[Interpolation coefficients]
The linearization \eqref{linear} can be generated from evaluations of $f_i$ as in \eqref{summation} at the Chebyshev nodes $\mu_{\ell}^*$ described in \eqref{interpolation-condition}. In this way, we do not explicitly require the functions $f_i$ in order to carry out the linearization.
\end{remark}

\section{Preliminaries for Chebyshev BiCG with exact preconditioning} \label{sec:Shifted-Systems-BiCG}

We consider a preconditioned short-term recurrence method where two Krylov subspaces are generated via a Lanczos biorthogonalization procedure as in the method BiCG. The biorthogonalization process requires the action of the system matrix, as well as its adjoint. Additionally, our setting uses the shift-and-invert preconditioner $(K-\sigma M)^{-1}$, for $\sigma \in (-a,a)$, and the corresponding adjoint preconditioner $(K-\sigma M)^{-T}$. Some preliminaries are described here as preparation, and the strategy for the application of the preconditioners is shown in Section~\ref{sec:Preconditioning}. The proposed method, Algorithm~\ref{alg:MCBiCG}, is derived and presented in Section~\ref{sec:Deriv-MCBiCG}, and numerical simulations follow in Section~\ref{sec:Helmholtz} and Section~\ref{sec:delay-problem}. 

Our right preconditioned system is of the form
\begin{subequations} \label{shifted}
\begin{eqnarray}
	&&(K - \mu M) (K - \sigma M)^{-1} \tilde{u}(\mu) = \tilde{b} \\ &\iff& 
	(K - \mu M + \sigma M - \sigma M) (K - \sigma M)^{-1} \tilde{u}(\mu) = \tilde{b} \\&\iff&
	\left(I + (-\mu+\sigma) M (K-\sigma M)^{-1} \right) \tilde{u}(\mu) = \tilde{b}, \label{shifted-c}
\end{eqnarray}
\end{subequations}
where $\tilde{u}(\mu) = (K-\sigma M)u(\mu)$ and 
\begin{align} \label{conv-reasoning}
	(K- \mu M)(K - \sigma M)^{-1} \approx I,  
\end{align}
for $\mu \approx \sigma$. The parameter $\sigma$ in the preconditioners is chosen based on the values of $\mu$ we are interested in and, thus, can be seen as a target parameter. Specifically, we chose $\sigma$ such that we approximate \eqref{our-problem} for many different values of $\mu$ in a neighborhood of $\sigma$. 

Krylov methods have been developed to approximate the solution to shifted linear systems of the same form as \eqref{shifted-c} in many prior works. See, for example, \cite{Freund:NLAproc92,Baumann2015NestedKM,FrommerGlassner98}, as well as \cite{Bahkos:2017:MULTIPREC}, where multiple shift-and-invert preconditioners were incorporated to build a rich search space in a GMRES framework, and \cite{ParksEtAl2006,SOODHALTER2014105}, where Krylov recycling techniques were utilized to solve shifted linear systems. Additionally, in \cite{BAUMANN2019316}, preconditioned Krylov subspace methods were considered for the time-harmonic elastic wave equation. Specifically, the solution to an equivalent linearized system was approximated, and shift-and-invert preconditioners with complex linear shifts were considered. 

The formulation \eqref{shifted-c} allows us to take advantage of the shift- and scaling-invariance properties of Krylov subspaces. More concretely, on the $j$th iteration of our proposed BiCG method Algorithm~\ref{alg:MCBiCG}, we seek approximations to $\tilde{u}(\mu)$ in~\eqref{shifted} from the Krylov subspace defined by
\begin{align} \label{krylov-subspace}
	\mathcal{K}_{j} \coloneqq \text{span} \left\{ \tilde{b}, M(K-\sigma M)^{-1} \tilde{b}, \ldots, \left(M(K-\sigma M)^{-1}\right)^{j-1} \tilde{b} \right\}.
\end{align}
Here $\mathcal{K}_j = \mathcal{\tilde{K}}_j$, where $\mathcal{K}_j = \mathcal{K}_j( M(K-\sigma M)^{-1}, \tilde{b} )$ and $\mathcal{\tilde{K}}_j$ is the Krylov subspace of dimension $j$ generated from the system matrix in~\eqref{shifted-c} and the vector $\tilde{b}$. In this way, approximating the solution to \eqref{shifted} for many values of $\mu$ in a BiCG setting requires just one basis matrix for $\mathcal{K}_j$ and, analogously, one basis matrix for the Krylov subspace $\mathcal{L}_j$, defined by
\begin{align} \label{adjoint-krylov-sub}
	\mathcal{L}_j \coloneqq \text{span} \left\{ \tilde{c}, \left( M(K-\sigma M)^{-1} \right)^T \tilde{c}, \ldots, \left(\left(M(K-\sigma M)^{-1} \right)^T \right)^{j-1} \tilde{c}\right\},
\end{align}
for $\tilde{b}^T \tilde{c} \neq 0$, $\tilde{c} \in \RR^{dn}$. Equivalently, we use one basis matrix for \eqref{krylov-subspace} to approximate solutions to \eqref{our-problem} for many values of $\mu$, as described in Section~\ref{sec:linearization}. To our knowledge, this is the first time a Chebyshev linearization has been combined with a Krylov subspace method to obtain approximate solutions to parameterized linear systems. 

In particular, after $j$ iterations, the Lanczos biorthogonalization generates matrices $V_j$, $W_j \in \RR^{dn \times j}$, $T_j \in \RR^{j \times j}$, and $\underbar{$T$}_j$, $\bar{T}^T_j \in \mathbb{R}^{(j+1) \times j}$ such that the relations
\begin{subequations} \label{relations}
\begin{alignat}{5} 
	M(K-\sigma M)^{-1} V_{j} &= V_{j} &&T_j &&+ \beta_j v_{j+1} e_j^T &&= V_{j+1} &&\underbar{$T$}_j, \label{relations-a} \\
	\big(M(K-\sigma M)^{-1}\big)^T W_{j} &= W_{j} &&T_j^T &&+ \gamma_j w_{j+1} e_j^T &&= W_{j+1} &&\bar{T}^T_j
\end{alignat}
\end{subequations}
hold, where the columns of $V_{j}$ span the subspace \eqref{krylov-subspace}. Analogously, the columns of $W_{j}$ span the subspace \eqref{adjoint-krylov-sub}, and the biorthogonalization procedure gives the relation
\begin{align}\label{iden}
	W_{j}^T V_{j} = I_{j},
\end{align}
where $I_j \in \mathbb{R}^{j \times j}$ is the identity matrix of dimension $j \times j$ and $e_j$ is the $j$th column of $I_j$.  Here the square matrix $T_j$ has the form
\begin{align} \label{square-T}
	T_{j} \coloneqq 
	\begin{bmatrix}
	\alpha_1 & \gamma_1 & & & \\
	\beta_1 & \ddots & & & \\
	& \ddots & & & \gamma_{j-1}\\\
	& & & \beta_{j-1}& \alpha_{j} \\ 
	\end{bmatrix} \in \mathbb{R}^{j \times j},
\end{align}
and the tridiagonal Hessenberg matrices $\underline{T}_j$ and $\bar{T}^T_j$ are given by
\begin{align} \label{Ts}
	\underline{T}_j \coloneqq 
	\begin{bmatrix}
	\alpha_1 & \gamma_1 & & \\
	\beta_1 & \ddots & \ddots& \\
	& \ddots &\ddots & \gamma_{j-1}\\\
	& & \beta_{j-1}& \alpha_j \\ 
	& & & \beta_j 
	\end{bmatrix}, \quad
	\bar{T}^T_j\coloneqq 
	\begin{bmatrix}
	\alpha_1 & \beta_1 & &  \\
	\gamma_1 & \ddots & \ddots&  \\
	& \ddots &\ddots & \beta_{j-1} \\
	& & \gamma_{j-1} & \alpha_j \\
	& & & \gamma_j 
	\end{bmatrix}.
\end{align}
Note, only the $j \times j$ principal submatrices of $\underline{T}_j$ and $\bar{T}^T_j$ are the transpose of each other. An advantage of the Lanczos biorthogonalization procedure is the so-called short-term recurrence of the Krylov basis vectors, i.e., that the matrices $\underline{T}_j$ and $\bar{T}^T_j$ in \eqref{relations} are tridiagonal. In this way, the basis vectors are computed recursively at each iteration of the algorithm. This choice of method is different from those in the previous works \cite{CorrentyEtAl,JarlebringCorrenty1} for parameterized systems based on companion linearization. The resulting linear systems in these works were solved in a GMRES setting, and they used a Gram-Schmidt orthogonalization process, i.e., a long-term recurrence.

The residual of the $i$th iterate generated from a BiCG procedure applied to the shifted preconditioned linear system \eqref{shifted}, denoted $\tilde{r}_i$, is orthogonal to the subspace $\mathcal{L}_i$ as in \eqref{adjoint-krylov-sub}, and $\tilde{r}_i$ is an element of the Krylov subspace $\mathcal{K}_{i+1}$ as in \eqref{krylov-subspace}. As established in~\cite{Frommer2003res}, $\tilde{r}_i$ is colinear to the residual vector of the $i$th iterate resulting from BiCG applied to the so-called \textit{seed} system given by
\begin{align} \label{seed}
	M(K-\sigma M)^{-1} u^{sd} = \tilde{b},
\end{align}
denoted $r_i$. This is described entirely in Section~\ref{sec:Deriv-MCBiCG} and used to derive Algorithm~\ref{alg:MCBiCG}. Analogously, we denote the adjoint seed system by
\begin{align} \label{adjoint-seed}
	(K-\sigma M)^{-T}M^T w^{sd} = \tilde{c}. 
\end{align}

Approximates to \eqref{our-problem} corresponding to each $\mu$ in an interval can be obtained with little extra computation using the {\em exact algorithm} described in Section~\ref{sec:Deriv-MCBiCG}, and this method is based on the colinearity of $r_i$ and $\tilde{r}_i$. Here the accuracy of the approximations depends on the dimension $j$ of the Krylov subspace from which these approximates originate. Additionally, the basis matrix for the Krylov subspace does not need to be stored if the values of $\mu$ are determined in advance. Our method offers a computational advantage over solving each parameterized system individually when we are interested in the solution to many different parameterized systems.

\section{An efficient application of a shift-and-invert preconditioner} \label{sec:Preconditioning}
Utilizing preconditioning in the context of Krylov subspace methods can lead to methods which are more efficient overall. This strategy is only suitable when the action of the preconditioner is cheap to apply. Our proposed method incorporates well-established shift-and-invert preconditioning. As we consider a BiCG setting, we require an efficient application of the preconditioner and its adjoint, i.e., we consider applying the action of $(K-\sigma M)^{-1}$ as well as the action of $(K-\sigma M)^{-T}$, $\sigma \in (-a,a)$. 

These preconditioners are effective for solving shifted systems of the form \eqref{companion} when the target parameter $\sigma$ is chosen close to the values of $\mu$ of interest. This is due to the relation~\eqref{conv-reasoning}. The BiCG method tends to experience fast converge when applied to linear systems of this form. In the following, we show how the structure can be exploited such that the action of these preconditioners can be computed efficiently. 

Consider approximating the solution to the system \eqref{companion}, incorporating right preconditioning with shift-and-invert preconditioner $(K-\sigma M)^{-1}$. The resulting linear system is of the form in \eqref{shifted}, expressed equivalently as
\begin{align} \label{last-in-shifted}
 (-\mu+\sigma)\left( \frac{1}{(-\mu+\sigma)} I + M(K-\sigma M)^{-1} \right)  \tilde{u}(\mu) = \tilde{b}.
 \end{align}
Note, the formulation in \eqref{last-in-shifted} was chosen in order to match the notation in the Exact Algorithm~\ref{alg:MCBiCG}, presented in Section~\ref{sec:Deriv-MCBiCG}. Specifically, the scalar $1/(-\mu+\sigma)$ is the coefficient of the identity matrix. This work considers a BiCG setting, and the linear system \eqref{last-in-shifted} incorporates a shift with a scalar multiple of the identity matrix. Equivalent to \eqref{relations}, the shifted relations
\begin{subequations} \label{shifted-relations}
\begin{eqnarray}
	V_j + (-\mu + \sigma) M(K-\sigma M)^{-1} V_j &=& V_{j+1} (\underline{I}_j + (-\mu + \sigma) \underline{T}_j ), \label{shifted-relations-a} \\
	W_j + (-\mu + \sigma) \big(M(K-\sigma M)^{-1} \big)^T W_j &=& W_{j+1} (\underline{I}_j + (-\mu + \sigma) \bar{T}^T_j ),	
\end{eqnarray}
\end{subequations}	
hold, where the matrix $\underline{I}_j \in \mathbb{R}^{(j+1) \times j}$ is an identity matrix of dimension $j \times j$ with an extra row of zeros, i.e., $V_{j+1} \underline{I}_j = V_j$. Note, the matrix $V_j$ is also a basis for the Krylov subspace generated by the matrix in \eqref{last-in-shifted} and the vector $\tilde{b}$ by the shift- and scaling-invariance properties of Krylov subspaces described in Section~\ref{sec:Shifted-Systems-BiCG}. 

Approximating the solution to \eqref{last-in-shifted} in our setting requires a basis for the Krylov subspaces~\eqref{krylov-subspace} and~\eqref{adjoint-krylov-sub}, obtained from multiplication with the matrices in~\eqref{seed} and~\eqref{adjoint-seed} at each iteration of the Lanczos biorthogonalization. The action of matrices $M$ and $M^T$ require one matrix-vector multiplication of size $n \times n$ for each product. Additionally, an efficient application of the preconditioner $(K-\sigma M)^{-1}$ can be performed via a block LU decomposition of the matrix $(K-\sigma M)\Pi$, where $\Pi$ is a permutation matrix, as described in~\cite{Amiraslani2009,KressRoman}. This process leads directly to an efficient application of $(K-\sigma M)^{-T}$. The procedure is described as follows. 

Let a permutation of the block columns of the matrix $(K-\sigma M)$ in \eqref{linear} be given by
\begin{align} \label{permutation-matrix}
	(K-\sigma M) \Pi = 
	\begin{bmatrix}
	I & & & & & -\frac{\sigma}{a}I\\
	-\frac{2\sigma}{a}I&I & & & &I \\
	I &-\frac{2\sigma}{a}I& I& & & \\
	& &\ddots & & & \\
	& & I& -\frac{2\sigma}{a}I& I& \\
	P_1& \cdots & P_{d-3}& (-P_d+P_{d+2})& (P_{d-1}+\frac{2\sigma}{a}P_d)& P_0\\
	\end{bmatrix},
\end{align}
where $\Pi\coloneqq \footnotesize \begin{bmatrix} & I_n \\ I_{(d-1)n}\end{bmatrix} \in \mathbb{R}^{dn \times dn}$ is an orthogonal matrix, and let $L_{\sigma} U_{\sigma} = (K-\sigma M) \Pi$ be a block LU decomposition, where
\[
	L_{\sigma} \coloneqq 
	\begin{bmatrix}
	I& & & & & \\
	-\frac{2\sigma}{a}I& I& & & & \\
	I&-\frac{2\sigma}{a}I &I & & & \\
	& & \ddots& & & \\
	& &I & -\frac{2\sigma}{a}I& I& \\
	P_1& \cdots & P_{d-3}& (-P_d+P_{d+2})& (P_{d-1}+\frac{2\sigma}{a}P_d)& P(\sigma)\\
	\end{bmatrix}
\]
and
\[
	U_{\sigma} \coloneqq
	\begin{bmatrix}
	I& & & & & -\tau_1(\sigma)I
	\\
	& I& & & & -\tau_2 (\sigma)I\\
	& & I& & & -\tau_3 (\sigma)I\\
	& & & \ddots& & \vdots\\
	& & & & I& -\tau_{d-1} (\sigma)I\\
	& & & & & I\\
	\end{bmatrix}.
\]

The matrix $U_{\sigma}^{-1}$ is identical to $U_{\sigma}$, except for a sign change in the first $d-1$ blocks in the last block column. Applying $L_{\sigma}^{-1}$ to a vector amounts to recursively calculating the first $d-1$ block elements, in addition to one linear solve with system matrix $P(\sigma) \in \mathbb{R}^{n \times n}$ when computing the last block row. This process of applying $L_{\sigma}^{-1}$ is equivalent to applying Gaussian elimination on a lower block triangular system. 

Thus, the action of the preconditioner $(K-\sigma M)^{-1}$ applied to a vector $y \in \RR^{dn}$ is given by
\begin{align} \label{s-and-i}
	(K-\sigma M)^{-1} y = \Pi U_{\sigma}^{-1} L_{\sigma}^{-1} y;
\end{align}
i.e., the preconditioner in \eqref{shifted} can be applied without computing or storing the large matrices $K$ and $M$ or $L_{\sigma}^{-1}$ and $U_{\sigma}^{-1}$. Note, the action of $P(\sigma)^{-1}$ in the application of $L_{\sigma}^{-1}$ can be done, for example, via one LU decomposition of $P(\sigma)$ performed before the start of the algorithm, as we consider in Algorithm~\ref{alg:MCBiCG} in Section~\ref{sec:Deriv-MCBiCG}. Alternatively we can apply the action of $P(\sigma)^{-1}$ approximately via an iterative method, as considered in Algorithm~\ref{alg:FMCBiCG} and presented in Section~\ref{sec:Flex}. 

The action of the adjoint preconditioner $(K-\sigma M)^{-T}$ can be applied in an analogous way, i.e.,
\begin{align} \label{s-and-i-T}
	(K-\sigma M)^{-T} y = L_{\sigma}^{-T} U_{\sigma}^{-T} \Pi ^{T} y,
\end{align}
which does not require the storage of matrices $K$ and $M$ or the triangular matrices $L_{\sigma}^{-T}$ and $U_{\sigma}^{-T}$. Additionally, the LU factorization of $P(\sigma)$ required in Algorithm~\ref{alg:MCBiCG} can be reused in the application of the adjoint preconditioner. Thus, the shift-and-invert preconditioner $(K-\sigma M)^{-1}$ is suitable in a BiCG setting.  

\section{Derivation of Preconditioned Chebyshev BiCG for parameterized linear systems} \label{sec:Deriv-MCBiCG}
In \cite{Frommer2003res} a BiCG algorithm was derived in order to solve a seed system of the form \eqref{seed}, as well as a shifted system of the form \eqref{last-in-shifted}, without requiring additional matrix-vector products for the iterates of the shifted system. Additionally, in \cite{PreCondAhmed}, a multishift BiCG algorithm with polynomial preconditioning was proposed to approximate the solutions for a family of shifted systems simultaneously. Such approaches are far less costly than solving each system of interest individually in a BiCG setting, without degrading convergence in general. 

We summarize below a derivation of the original algorithm from \cite{Frommer2003res}, adapted to our shifted system stemming from a companion linearization formed from a Chebyshev interpolation as in \eqref{companion}. The method is based on the observation that the residual vectors of the seed system can be used to generate a basis for the Krylov subspace~\eqref{krylov-subspace}. The description which follows serves to clarify the steps in Algorithm~\ref{alg:MCBiCG}, which we will refer to as the Exact Algorithm~\ref{alg:MCBiCG} or the {\em exact algorithm}. 

Lines $1-8$ of the Exact Algorithm~\ref{alg:MCBiCG} correspond to a variant of the standard BiCG method. This implementation is based on a Lanczos biorthogonalization, applied to the seed system~\eqref{seed}. The coupled two-term recurrence formulation used here is based on an implicitly formed LU decomposition of the tridiagonal matrix $T_j$ as in \eqref{relations}; see \cite{IteratMSparseSys} for a detailed description. The search direction vectors in this formulation are updated as
\begin{align} \label{search-update}
	v_{i+1}^* = r_{i} - \beta_{i} v_{i}^*, \qquad w_{i+1}^* = s_{i} - \bar{\beta}_i w_{i}^*,  
\end{align}
where the residual vectors
\begin{align} \label{res-vecs}
r_i \coloneqq \tilde{b} - M(K-\sigma M)^{-1} u^{sd}_i, \qquad s_i \coloneqq \tilde{c} - (K-\sigma M)^{-T}M^T w^{sd}_i
\end{align}
are updated recursively (see line 7 of the algorithm), and $u_i^{sd}, w_i^{sd} \in \RR^{dn}$ are approximations to the seed system and adjoint seed systems given in \eqref{seed} and \eqref{adjoint-seed}, respectively, with $u_0^{sd} \coloneqq 0$, $w_0^{sd} \coloneqq 0$. The search directions have the property $(w_i^*)^T M (K-\sigma M)^{-1} v_j^*=0$, for $i \neq j$, and residual vectors are orthogonal, i.e., $s_i^T r_j = 0$, for $i\neq j$. The approximate solution to the seed system is updated in line 8 as 
\begin{align} \label{sol-update}
	u_{i+1}^{sd} = u_{i}^{sd} + \alpha_i v_{i+1}^*.
\end{align}

Let $\tilde{u}_i (\mu) \in \RR^{dn}$ be the $i$th approximate solution to the shifted preconditioned linear system \eqref{last-in-shifted}, obtained from the BiCG method. Since $u_{i}^{sd}, \tilde{u}_{i}(\mu)$ are elements of the Krylov subspace $\mathcal{K}_{i}$ as in \eqref{krylov-subspace}, we can express these approximate solutions as
\[
	u_{i}^{sd} = p_{i-1} \left( M(K-\sigma M)^{-1} \right) \tilde{b}, \qquad
	\tilde{u}_{i}(\mu) = \tilde{p}_{i-1} \left(\frac{1}{(-\mu+\sigma)} I + M (K-\sigma M)^{-1} \right)\tilde{b},
\]
for $i =1,2,\ldots$, where $p_{i-1}$, $\tilde{p}_{i-1}$ are polynomials of degree less than or equal to $i-1$. Similarly, the residual $r_i$ in \eqref{res-vecs} and the residual of the $i$th iterate obtained from BiCG applied to the shifted system, denoted $\tilde{r}_i$, can be expressed as
\begin{align} \label{r-seed}
	r_i = q_i \left( M(K-\sigma M)^{-1} \right) \tilde{b}, \qquad \tilde{r}_i &= \tilde{q}_i \left( \frac{1}{(-\mu+\sigma)} I + M (K-\sigma M)^{-1} \right) \tilde{b},
\end{align}
where $q_i(t) \coloneqq 1 - t p_{i-1} (t)$, $\tilde{q}_i (t) \coloneqq 1 - t \tilde{p}_{i-1}(t)$, and 
\begin{align} \label{zero-eval}
	q_i(0) = \tilde{q}_i(0) = 1. 
\end{align}
Here the polynomials $q_i$ and $\tilde{q}_i$ are of degree less than or equal to $i$. The following lemma characterizes the relation between $r_i$ and $\tilde{r}_i$, originally established in \cite{Frommer2003res}. 

\begin{lemma} \label{lemma1}
Let $r_i$ and $\tilde{r}_i$ as in \eqref{r-seed} be the residuals at iteration $i$ resulting from the BiCG method applied to \eqref{seed} and \eqref{last-in-shifted}, respectively. Then, there exists $\zeta_i \in \mathbb{R} $ such that $r_i = \zeta_i \tilde{r}_i$, i.e., $r_i$ and $\tilde{r}_i$ are colinear for $i=1,\ldots,j$. 
\end{lemma}

We omit the proof of Lemma~\ref{lemma1} as it is analogous to the proof of \cite[Theorem 1]{Frommer2003res}. Briefly, the result is shown by noting the shift-invariance of Krylov subspaces as well as the property of the residual of BiCG iterates: $r_i, \tilde{r}_i \in (\mathcal{L}_i)^{\perp} \cap \mathcal{K}_{i+1}\eqqcolon  \mathcal{J}_i$, where $\mathcal{L}_i$ is as in \eqref{adjoint-krylov-sub} and $\mathcal{K}_{i+1}$ is as in \eqref{krylov-subspace}. The biorthogonality condition \eqref{iden} implies that the dimension of $\mathcal{J}_i$ is 1, and the result follows. 
 
By Lemma \ref{lemma1}, the relation between the residual of the shifted system and the seed system is given by 
\begin{align} \label{plus}
	\tilde{r}_i = \frac{1}{\zeta_i} r_i,
\end{align}
or, equivalently, with residual polynomials,
\begin{align} \label{r-tilde-rel}
	\tilde{r}_i = \tilde{q}_i \left( \frac{1}{-\mu + \sigma} I + M (K-\sigma M)^{-1} \right) \tilde{b} = \frac{1}{\zeta_i} q_i \left( M (K-\sigma M)^{-1} \right) \tilde{b}.
\end{align}
The equality in \eqref{r-tilde-rel} can be expressed as a function of $t$, i.e., $\tilde{q}_i \left( 1/(-\mu+\sigma)+ t \right) = (1/\zeta_i) q_i (t),$ and, paired with the equality in \eqref{zero-eval}, gives
\begin{align} \label{triangle}
	\zeta_i = q_i \left( \frac{-1}{-\mu + \sigma} \right). 
\end{align}
Note, the above gives an exact expression for the colinearity coefficient $\zeta_i$ in Lemma \ref{lemma1}, completely determined from the residual of the seed system. Thus, from \eqref{plus} and \eqref{triangle}, we can express the $i$th residual resulting from BiCG applied to the shifted system from the residual of the $i$th iterate obtained from BiCG applied to the seed system. In other words, we can obtain the residual vectors corresponding to many different shifted systems from one execution of the algorithm and, as a result, update the search vectors in \eqref{search-update} to approximate the solution to many shifted systems as in \eqref{sol-update}. 

The derivation which follows serves to clarify lines $9-14$ of the Exact Algorithm~\ref{alg:MCBiCG}, where our algorithm is applied to \eqref{last-in-shifted} for a set of shifts $\{ \mu_l \}$, $l=1,\ldots,k$. In the algorithm, we denote the particular $\zeta_i$ for each shift $\mu_l$ as $\zeta_i (\mu_l)$, but use the notation $\zeta_i\coloneqq\zeta_i(\mu)$ here for simplicity. 

The residual of the BiCG iterates applied to the seed system satisfy the following three-term recursion:
\begin{align} \label{star}
	r_{i+1} = -\alpha_i M(K-\sigma M)^{-1} r_i + \frac{\beta_i \alpha_i}{\alpha_{i-1}} r_{i-1} + \left( 1- \frac{\beta_i \alpha_i}{\alpha_{i-1}} \right) r_i,
\end{align} 
where we have inserted the update formula $v_{i+1}^* = r_i - \beta_i v_i^*$ from line 3 of the algorithm into the computation in line 7, i.e., 
\begin{align*}
	r_i &= r_{i-1} - \alpha_{i-1} M (K-\sigma M)^{-1} v_i^* \\
	&= r_{i-1} - \alpha_{i-1} M(K-\sigma M)^{-1} \left( \frac{1}{\beta_i} (r_i - v_{i+1}^*) \right) \\
	&= r_{i-1} - \frac{\alpha_{i-1}}{\beta_i} M (K-\sigma M)^{-1} r_i + \frac{\alpha_{i-1}}{\beta_i \alpha_i} (r_i - r_{i+1}),
\end{align*}
and used the recursive update formula $M(K-\sigma M)^{-1} v_{i+1}^* = (1/\alpha_i ) (r_{i} - r_{i+1})$ in the last equality. The relation \eqref{star} can be expressed with the residual polynomial from \eqref{r-seed} as
\begin{align} \label{poly-rel}
	q_{i+1} (t) = -\alpha_i t q_i(t) + \frac{\beta_i \alpha_i}{\alpha_{i-1}} q_{i-1}(t) + \left( 1 - \frac{\beta_i \alpha_i}{\alpha_{i-1}} \right) q_i(t).
\end{align}
Specifically, taking $t = -1/(-\mu+\sigma)$ in \eqref{poly-rel} and using the equality in \eqref{triangle}, gives 
\begin{align} \label{zeta-recurrence}
	\zeta_{i+1} = \left( 1 - \alpha_i \Big(\frac{-1}{-\mu+\sigma}\Big) - \frac{\beta_i \alpha_i}{\alpha_{i-1}}\right) \zeta_i + \frac{\beta_i \alpha_i}{\alpha_{i-1}} \zeta_{i-1}, 
\end{align}
i.e., a recurrence for the colinearity coefficients in Lemma~\ref{lemma1}, incorporated in line 10 of the algorithm. Thus, from \eqref{plus}, \eqref{star}, and \eqref{zeta-recurrence}, the three-term recurrence of the residual vectors resulting from the BiCG method applied to the shifted system is given by
\begin{align*}
	\tilde{r}_{i+1} &= \frac{1}{\zeta_{i+1}} \left( -\zeta_i \alpha_i M (K-\sigma M)^{-1} \tilde{r}_i + \zeta_{i-1} \Big(\frac{\beta_i \alpha_i}{\alpha_{i-1}}\Big) \tilde{r}_{i-1} + \zeta_i \Big( 1 - \frac{ \beta_i \alpha_i}{\alpha_{i-1}} \Big) \tilde{r}_i  \right) \\
 &= -\frac{\zeta_i \alpha_i}{\zeta_{i+1}} \left( \frac{1}{-\mu+\sigma} I + M (K-\sigma M)^{-1} \right) \tilde{r}_i 
	+ \frac{\zeta_{i-1}}{\zeta_{i+1}} \left(\frac{\beta_i \alpha_i}{\alpha_{i-1}}\right) \tilde{r}_{i-1} 
	+ \left(  1 - \frac{\zeta_{i-1}}{\zeta_{i+1}} \Big( \frac{\beta_i \alpha_i}{\alpha_{i-1}} \Big) \right) \tilde{r}_i,  	
\end{align*}
where we have used the relation \eqref{zeta-recurrence} to obtain the third term of the summation in the last equality. Equivalently, we can express the $(i+1)$st residual generated from approximating the solution to \eqref{last-in-shifted} with BiCG with the recurrence 
\[
	\tilde{r}_{i+1} = -\tilde{\alpha}_i \left( \frac{1}{-\mu+\sigma} I + M (K-\sigma M)^{-1} \right) \tilde{r}_i + \frac{\tilde{\beta}_i \tilde{\alpha}_i}{\tilde{\alpha}_{i-1}} \tilde{r}_{i-1} + \left( 1 - \frac{\tilde{\beta}_i \tilde{\alpha}_i}{\tilde{\alpha}_{i-1}} \right) \tilde{r}_i,
\]
where the coefficients $\tilde{\alpha}_i$ and $\tilde{\beta}_i$ are defined as
\begin{align} \label{shift-coeffs}
\tilde{\alpha}_i \coloneqq -\alpha_i \left( \frac{\zeta_i}{\zeta_{i+1}} \right), \qquad
\tilde{\beta}_i \coloneqq \left( \frac{\alpha_i}{\tilde{\alpha}_i} \right) \left( \frac{\tilde{\alpha}_{i-1}}{\alpha_{i-1}} \right) \frac{\zeta_{i-1}}{\zeta_{i+1}} \beta_i = \left( \frac{\zeta_{i-1}}{\zeta_i} \right)^2 \beta_i
\end{align}
and updated in line 11 of the algorithm. Note, initializing with parameters $\zeta_0 = \zeta_1 = 1$ in~\eqref{zeta-recurrence} ensures that the formulation described above holds for the corresponding seed system (cf.~\eqref{star}). Analogous to \eqref{search-update}, we compute the search vectors for solving the shifted systems as 
\begin{align} \label{shifted-update-vec}
	\tilde{v}_{i+1} = \frac{1}{\zeta_i} r_i - \tilde{\beta}_i \tilde{v}_i, 
\end{align}
in line 12 and update the approximation to each shifted preconditioned system in line 13 as 
\begin{align} \label{tilde-u-update}
\tilde{u}_{i+1}(\mu) = \tilde{u}_{i}(\mu) + \tilde{\alpha}_i \tilde{v}_{i+1},
\end{align}
i.e., the shifted equivalent of the update described in \eqref{sol-update}. Lines 15-25 of the algorithm ensure that the approximations $\tilde{x}_j(\mu_l)$ to the linear system \eqref{our-problem} from the Krylov subspace of dimension $j$ have relative residual norm below a certain tolerance $tol$, for $l=1,\ldots,k$.

The Exact Algorithm~\ref{alg:MCBiCG} applies the shift-and-invert preconditioners $(K-\sigma M)^{-1}$ and $(K-\sigma M)^{-T}$ via a block LU decomposition as described in Section~\ref{sec:Preconditioning}. Note, for each update of the solution to the seed system, only some additional scalar operations and vector additions are required to update the approximations to $\tilde{u}(\mu_l)$ as in \eqref{last-in-shifted} for each $\mu_l$. This is due to the colinearity of the residuals $r_i$ and $\tilde{r}_i$, as described in Lemma~\ref{lemma1}. Furthermore, the Exact Algorithm~\ref{alg:MCBiCG} does not require the storage of the residual vectors in~\eqref{res-vecs} at each iteration, as long as the values of $\mu_l$, $l=1,\ldots,k$ are determined before the algorithm is executed. This allows for a method with low memory consumption, even when the degree $d$ of the Chebyshev approximation is large. If $\rho_i$ in line 2 vanishes, the algorithm has a breakdown. This scenario never occurred in our experiments.   

\begin{algorithm}[p!]
\SetAlgoLined
\SetKwInput{Input}{Input}\SetKwInput{Output}{Output}\SetKwInput{Initialize}{Initialize}
\SetKwFor{FOR}{for}{do}{end for}
\SetKwFor{IF}{if}{then}{end if}
\SetKwFor{FUN}{function}{}{end function}
\SetKwFor{WHILE}{while}{do}{end while}
\Input{$P_{\ell}$, $\ell = 0,1,\ldots,d$, as in \eqref{chebyshev} (Chebyshev coefficients) \\ $L_P U_P$, decomposition of $P(\sigma) \in \mathbb{R}^{n \times n}$ for $L_{\sigma}^{-1}$ \\$\tilde{b} \in \RR^{dn}$ as in \eqref{b-tilde}, $\tilde{c} \in \RR^{dn}$ such that $\tilde{b}^T \tilde{c} \neq 0$ \\ $\sigma \in \RR$ as in \eqref{shifted} (target parameter), $tol$\\ $\{ \mu_l \}_{l=1,\ldots,k}$, $\mu_l \in \mathbb{R}$, $|\mu_1-\sigma| \leq |\mu_2-\sigma| \leq \cdots \leq |\mu_k-\sigma|$ (ordered set of shifts) \\ $A(\mu_l)$, $l=1,\ldots,k$, as in \eqref{summation}, $\zeta_{-1}(\mu_l) = \zeta_{0}(\mu_l) = 1$, $l=1,\ldots,k$}
\Output{Approx. sol. $\tilde{x}_{j}(\mu_l)$, $l=1,\ldots,k$, to \eqref{our-problem}, from subspace of dim. $j$}
\Initialize{$\rho_{-1}=1$, $\alpha_{-1} = 1$, $\omega (\mu_l)=-1/(-\mu_l + \sigma)$, $l=1,\ldots,k$ \\
$v_{0}^* = w_{0}^* = 0 \in \RR^{dn}$ (search direction, seed system)
\\ $\tilde{v}_{0}(\mu_l)=0 \in \RR^{dn}$, $l=1,\ldots,k$ (search directions, shifted systems) \\
$u_0^{sd} = 0 \in \RR^{dn}$ (approx. to seed system) \\ $\tilde{u}_{0}(\mu_l)=0 \in \RR^{dn}$, $l=1,\ldots,k$ (approx. to shifted systems) \\
$r_{0} = \tilde{b}$, $s_{0} = \tilde{c}$ (residual vectors)} 
\FOR{$i = 0,1,2, \dots$ }{
$\rho_{i}=(r_{i})^T s_{i}$, $\beta_{i}=-\rho_i/\rho_{i-1}$ \\
$v_{i+1}^* = r_{i} - \beta_{i} v_{i}^*$, $w_{i+1}^* = s_{i} - \bar{\beta}_{i} w_{i}^*$ \\
Compute $\hat{v}_1$ such that $\hat{v}_1 =  M \left((K-\sigma M)^{-1} v_{i+1}^* \right)$ as in \eqref{s-and-i} \\
$\alpha_{i} = \rho_i / ((w_{i+1}^*)^T \hat{v}_1)$ \\
Compute $\hat{v}_2$ such that $\hat{v}_2 = (K-\sigma M)^{-T} (M^T w_{i+1}^*)$ as in \eqref{s-and-i-T} \\
$r_{i+1} = r_{i} - \alpha_{i} \hat{v}_1$, $s_{i+1} = s_{i} - \bar{\alpha}_{i} \hat{v}_2$ \\
$u_{i+1}^{sd} = u_{i}^{sd} + \alpha_{i} v_{i+1}^*$ (approx. sol. to seed system from subspace of dim $i+1$)\\
\FOR{l=1,\ldots,k} 
{$\zeta_{i+1}(\mu_l) = (1 - \alpha_i \omega(\mu_l) - \frac{\beta_i \alpha_i}{\alpha_{i-1}}) \zeta_i(\mu_l)+ \frac{\beta_i \alpha_i}{\alpha_{i-1}} \zeta_{i-1}(\mu_l)$ as in \eqref{zeta-recurrence} \\
$\tilde{\alpha}_i(\mu_l) = -\alpha_i \left(\frac{\zeta_i(\mu_l)} {\zeta_{i+1}(\mu_l)}\right)$, $\tilde{\beta}_i(\mu_l) = \left(\frac{\zeta_{i-1}(\mu_l)}{\zeta_i(\mu_l)}\right)^2 \beta_i$ as in \eqref{shift-coeffs}  \\
$\tilde{v}_{i+1}(\mu_l) = \left( \frac{1}{\zeta_i(\mu_l)} \right) r_i - \tilde{\beta}_i(\mu_l) \tilde{v}_{i}(\mu_l)$ as in \eqref{shifted-update-vec} \\ 
$\tilde{u}_{i+1}(\mu_l) = \tilde{u}_{i}(\mu_l) + \tilde{\alpha}_i(\mu_l) \tilde{v}_{i+1}(\mu_l)$ as in \eqref{tilde-u-update}}
Compute $\tilde{x}_{i+1}(\mu_k) = \textsc{PostProcess}(\omega(\mu_k),\tilde{u}_{i+1}(\mu_k))$ \\
Compute $res = \norm{A(\mu_k) \tilde{x}_{i+1}(\mu_k) - b}/\norm{b}$ \\
\IF {$res \leq tol$}
{Set $\tilde{x}_{j}(\mu_k)=\tilde{x}_{i+1}(\mu_k)$, $j = i+1$ \\
\FOR{$l=1,\ldots,k-1$}{$\tilde{x}_{j}(\mu_l)=$\textsc{ PostProcess}($\omega(\mu_l)$,$\tilde{u}_j (\mu_l)$) \\
 }
 \IF {$\norm{A(\mu_l) \tilde{x}_{j}(\mu_l) - b}/\norm{b} \leq tol$, $l=1,\ldots,k-1$} 
{return}}
}
\FUN {$\tilde{x}_{i+1}(\mu)=$\textsc{ PostProcess}($\omega(\mu)$,$\tilde{u}_{i+1}(\mu)$)} {
Compute $\hat{v}_3 = \omega(\mu) (K-\sigma M)^{-1} \tilde{u}_{i+1}(\mu)$ as in \eqref{s-and-i} \\
Set $\tilde{x}_{i+1}(\mu)= \hat{v}_3(1:n)$ \\
}
\caption{Preconditioned Chebyshev BiCG for parameterized linear systems}
\label{alg:MCBiCG}
\end{algorithm}

\begin{remark}[Adjoint parameterized system]
The Exact Algorithm~\ref{alg:MCBiCG} can be used to approximate the shifted right preconditioned adjoint linear system given by
\begin{subequations} \label{adjoint}
\begin{eqnarray}
	&& (-\mu + \sigma) \left( \frac{1}{(-\mu + \sigma)}I + (K-\sigma M)^{-T}M^T \right) \tilde{w} (\mu) = \tilde{c} \\
	&\iff& \left( (K-\sigma M)^T + (-\mu + \sigma) M^T \right) \tilde{w}(\mu) = (K-\sigma M)^T \tilde{c} \\
	&\iff& (K-\mu M)^T \tilde{w}(\mu) = (K-\sigma M)^T \tilde{c},
\end{eqnarray}
\end{subequations}
where $\tilde{w}(\mu) \in \RR^{dn}$. This system can be viewed as the shifted version of the adjoint seed system~\eqref{adjoint-seed}. Solutions to the adjoint system $A(\mu)^T z(\mu) = c$, $c \in \RR^n$ cannot be recovered from the above system, due to the structure of the adjoint of the linear system \eqref{linear}. Specifically, the solution vector in \eqref{adjoint} does not contain $z(\mu)$ since the Chebyshev interpolation coefficients of $A^T(\mu)$ appear in the last block column of the matrix $(K-\mu M)^T$, as opposed to the last block row of the matrix $(K-\mu M)$ in \eqref{linear}.
\end{remark} 

\section{Simulation of a parameterized Helmholtz equation} \label{sec:Helmholtz}
To highlight the capabilities of our method, we consider a Helmholtz equation, which describes the propagation of waves. Successful approaches for solving the Helmholtz equation have been considered in prior works such as \cite{BaylissEtAl1983,ElmanEtAl2001}, as well as in \cite{Erlangga08,ErlanggaEtAl2004}, where preconditioning was combined with fast iterative solvers. In particular, we consider the parameterized Helmholtz equation given by
\begin{subequations}\label{pde}
\begin{alignat}{2} 
	\left( \nabla^2 + \sin^2(\mu) \alpha(x) + \mu^2 + \cos^2(\mu) \beta(x) \right) u(x) &= h(x), \quad &&x \in \Omega, \\
	u(x) &= 0, &&x \in \partial \Omega,
\end{alignat}
\end{subequations}
where $\alpha(x)=1+\sin(x_1)$, $\beta(x) = 1+\cos(x_2)$, $h(x)=\exp(-x_1 x_2)$, and $\Omega \subset ([0,1] \times [0,1])$ is as in Figure~\ref{fig2}. The parameter $\mu$ in \eqref{pde} can be interpreted as a material parameter. It is of interest to approximate the solution $u(x)$ for a variety of different values $\mu$. 

\begin{figure}[h!]
\centering
\begin{subfigure}[T]{0.45\columnwidth}
	\centering
	\includegraphics[]{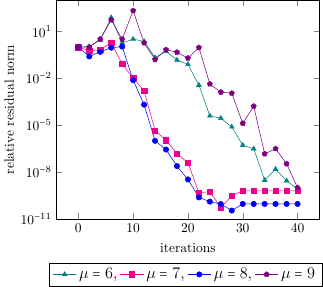}
\caption{Relative residual norm in terms of iterations}
\label{fig1a}
\end{subfigure}
\hspace{.5cm}
\begin{subfigure}[T]{0.45\columnwidth}
	\centering
	\includegraphics[]{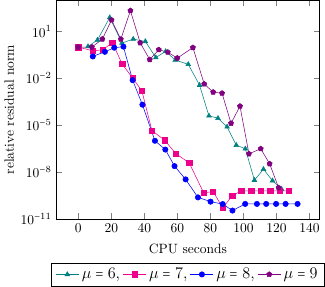}
\caption{Relative residual norm in terms of CPU sec}
\label{fig1b}
\end{subfigure}
\caption{Convergence for approximating the parameterized Helmholtz equation \eqref{pde}, generated from one execution of the Exact Algorithm~\ref{alg:MCBiCG} with $\sigma=7.5$ and evaluated for different values of $\mu \in [-10,10]$. Here $n=243997$, $d=50$, $tol=10^{-9}$, and relative residual norm \eqref{rel-res-pde}.}
\label{fig1}
\end{figure}

\begin{figure}[h!]
\centering
\includegraphics[]{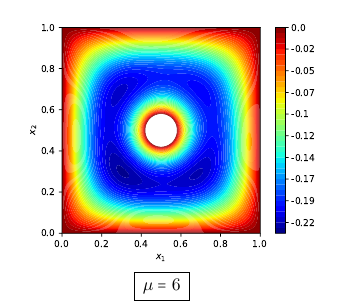}
\hspace{1cm}
\includegraphics[]{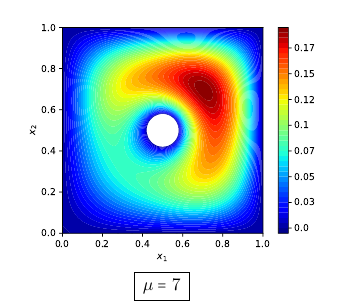}

\includegraphics[]{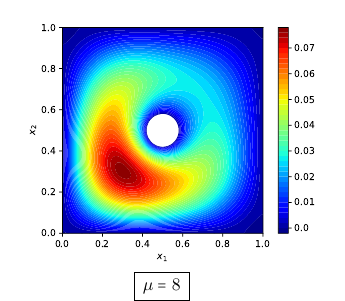}
\hspace{1cm}
\includegraphics[]{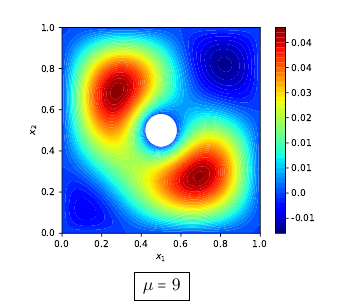}
\caption{Solutions of a Helmholtz equation corresponding to the simulation in Figure~\ref{fig1a} and one execution of the Exact Algorithm~\ref{alg:MCBiCG} with $\sigma = 7.5$.}
\label{fig2}
\end{figure}

\begin{figure}[h!]
\centering
\begin{subfigure}[T]{0.45\columnwidth}
	\centering
	\includegraphics[]{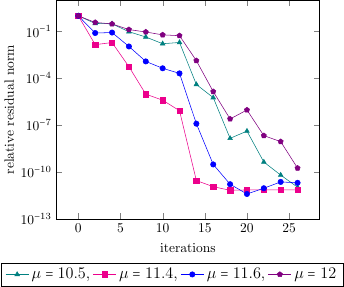}
	\caption{Relative residual norm in terms of iterations}
	\label{fig3a}
\end{subfigure}
\hspace{.5cm}
\begin{subfigure}[T]{0.45\columnwidth}
	\centering
	\includegraphics[]{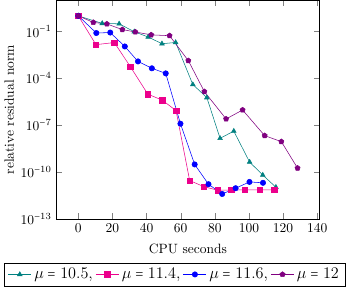}
	\caption{Relative residual norm in terms of CPU sec} 
	\label{fig3b}
\end{subfigure}
\caption{Convergence for approximating the parameterized Helmholtz equation \eqref{pde}, generated from one execution of the Exact Algorithm~\ref{alg:MCBiCG} with $\sigma=11.25$ and evaluated for different values of $\mu \in [-15,15]$. Here $n=243997$, $d=64$, $tol=10^{-9}$, and relative residual norm \eqref{rel-res-pde}.}
\label{fig3}
\end{figure}

\begin{figure}[h!]
\centering
\begin{subfigure}[T]{0.45\columnwidth}
	\centering
	\includegraphics[]{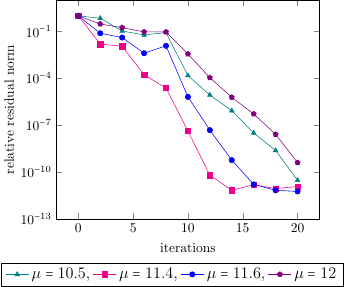}
	\caption{Relative residual norm in terms of iterations}
	\label{figg3ba}
\end{subfigure}
\hspace{.5cm}
\begin{subfigure}[T]{0.45\columnwidth}
	\centering
	\includegraphics[]{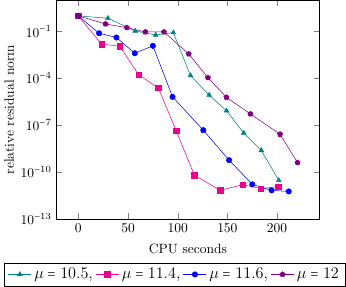}
	\caption{Relative residual norm in terms of CPU sec} 
	\label{figg3bb}
\end{subfigure}
\caption{Convergence for approximating the parameterized Helmholtz equation \eqref{pde}, generated from one execution of the Exact Algorithm~\ref{alg:MCBiCG} with $\sigma=11.25$ and evaluated for different values of $\mu \in [-40,40]$. Here $n=243997$, $d=124$, $tol=10^{-9}$, and relative residual norm \eqref{rel-res-pde}.}
\label{figg3b}
\end{figure}

\begin{figure}[h!]
\centering
\includegraphics[]{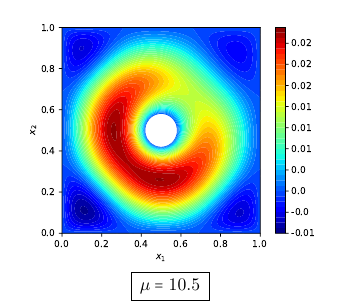}
\hspace{1cm}
\includegraphics[]{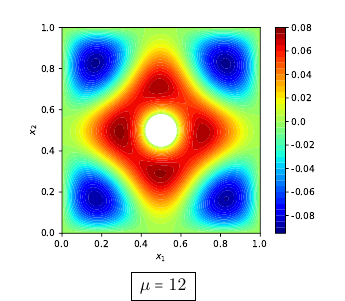}
\caption{Solutions of a Helmholtz equation corresponding to the simulation in Figure~\ref{fig3a} and one execution of the Exact Algorithm~\ref{alg:MCBiCG} with $\sigma = 11.25$.}
\label{fig4}
\end{figure}

Consider a discretization of \eqref{pde} which is of the same form as \eqref{our-problem}, i.e., 
\begin{align} \label{disc}
	A(\mu) \coloneqq A_0 + \sin^2(\mu) A_1 + \mu^2 A_2 + \cos^2(\mu) A_3,
\end{align}
where $A_0,\ldots,A_3$ arise from a finite element method (FEM) discretization and $b$ is the corresponding load vector.\footnote{The matrices and vector were generated using the finite element software FEniCS \cite{Alnaes:2015:FENICS}.} Approximating $A(\mu)$ \eqref{disc} with a Chebyshev interpolation leads to a parameterized linear system of the form \eqref{approx}, where $P(\mu) \approx A(\mu)$. We consider an approximate solution to the shifted preconditioned system \eqref{shifted}, based on a linearization of $P(\mu)$ as in \eqref{linear}. The resulting approximation to the companion linearization and, equivalently, to $x(\mu)$, is obtained for many values of the parameter $\mu$ via one execution of the Exact Algorithm~\ref{alg:MCBiCG}. The relative residual norm at iteration $i$,
\begin{align} \label{rel-res-pde}
	\frac{\norm{A(\mu) \tilde{x}_i (\mu) - b}}{\norm{b}}, 
\end{align}
is computed for a variety of $\mu$ with $A(\mu)$ \eqref{disc}. The results of this experiment are in Figures~\ref{fig1}-\ref{fig4}. Here, the nonlinear functions $\sin^2(\mu)$ and $\cos^2(\mu)$ are approximated using \texttt{Chebfun} in Matlab with truncation parameter $d$ as in \eqref{chebyshev}; see \cite{ChebfunGuide}. All simulations in this paper were carried out on a 2.3 GHz Dual-Core Intel Core i5 processor with 16 GB RAM. The software for all examples in this paper were implemented in Matlab, and we made them available online.\footnote{\texttt{https://github.com/siobhanie/ChebyshevBiCG}} 

Our experiments show that the preconditioned BiCG method leads to an accurate approximation of the linear system corresponding to a discretization of \eqref{pde} for values of the parameter $\mu$ close to the target parameter $\sigma$, and one execution of the Exact Algorithm~\ref{alg:MCBiCG} leads to a large variety of solution approximations. Specifically, as the algorithm is run to solve the seed system \eqref{seed}, each additional approximation corresponding to values of the parameter $\mu_l$, $l=1,\ldots,k$, is updated using only scalar operations and vector additions. Before each of the approximates to \eqref{our-problem} is returned, a final application of the preconditioner is performed. The basis for the Krylov subspace \eqref{krylov-subspace} does not need to be stored if the values of the parameter $\mu$ are determined in advance.

The simulation visualized in Figure~\ref{fig1a} gives access to the solutions corresponding to all values of $\mu \in [6,9]$, though we visualize just four of these. Analogously, the simulations shown in Figure~\ref{fig3a} and Figure~\ref{figg3ba} provide all solutions corresponding to $\mu \in [10.5,12]$. Note, the relative residual norms of the approximate solutions are all below a prescribed tolerance. 

In Figure~\ref{fig1b}, Figure~\ref{fig3b}, and Figure~\ref{figg3bb}, we see the benefit of the short-term recurrence property of the Exact Algorithm~\ref{alg:MCBiCG}, i.e., the roughly constant cost of each iteration. Specifically, we plot the relative residual norm in terms of CPU seconds, where we measure CPU time after the precomputation steps. This feature is especially useful for simulations which require many iterations until convergence and is not present in methods with a long-term recurrence, e.g., GMRES and the full orthogonal method (FOM). Note, though one execution of the Exact Algorithm~\ref{alg:MCBiCG} gives approximations for many different values of $\mu$, each convergence curve here corresponds to a separate run. 

From comparing the simulations in Figure~\ref{fig3} and Figure~\ref{figg3b}, we see that the cost of the approach is largely independent of the degree $d$ of the Chebyshev interpolation. More precisely, the interpolation corresponding to the simulation in Figure~\ref{fig3} was performed on the interval $[-15,15]$ with $d=64$, leading to a companion linearization of dimension $64 n \times 64 n$. Similarly, the simulation in Figure~\ref{figg3b} was performed on the interval $[-40,40]$ with $d=124$, where the companion linearization had dimension $124 n \times 124 n$. The experiments in Figure~\ref{fig3a} and Figure~\ref{figg3ba} converged in roughly the same number of iterations. The simulation in Figure~\ref{figg3bb} required approximately twice as many matrix-vector products with a matrix of dimension $n \times n$ as the one in Figure~\ref{fig3b} and took roughly twice as many CPU seconds as a result. Note, the cost of the application of the preconditioner is the same in both of these simulations as one LU decomposition of $P(\sigma) \in \mathbb{R}^{n \times n}$ is performed in the precomputation step. By performing the interpolation on a larger interval, we have access to a greater variety of solutions. However, only the solutions corresponding to values of $\mu$ close to the target $\sigma$ converge quickly. 

\begin{remark}[Magnitude of the parameter $\mu$]
The Chebyshev coefficients can be approximated efficiently using a discrete cosine transform. Thus, we can interpolate $A(\mu)$ on a large interval with little additional cost. While doing so leads to a larger linearization, it also allows for approximations to \eqref{our-problem} for large values of $\mu$. Considering values of $\mu$ close enough to $\sigma$ leads to convergence in $j$ iterations, where $j << dn$, with the additional benefit of a constant low cost per iteration. The experiments in this section were designed in order to show a wide variety of solutions. 
\end{remark}

\section{Simulation of a transfer function of a time-delay system} \label{sec:delay-problem}
Consider the solution to the linear system \eqref{our-problem}, where
\begin{align} \label{time-delay-sys}
	A(\mu) \coloneqq -\mu I + A_0 + A_1 e^{-\mu}
\end{align}
with random matrices $A_0$, $A_1 \in \mathbb{R}^{n \times n}$, random vector $b\in\mathbb{R}^n$, and $n=80$. The solution to this system is the transfer function of the time-delay system described by
\begin{align*}
	\dot{x}(t) &= A_0 x(t) + A_1 x (t - \tau) - b v(t), \\
	y(t) &= C^T x(t).
\end{align*}
Specifically, the transfer function is obtained by applying the Laplace transform to the state equation with $x(0)=0$. In this formulation, $\mu$ is the Laplace variable; see \cite{Gu:2003:STABILITY,Michiels2011,Michiels:2007:STABILITYBOOK}. The vector $b \in \RR^n$ is the external force, $x(t) \in \RR^n$ is the state vector, $v(t)$ is the input, $y(t)$ is the output and $\tau > 0$ is the delay. Without loss of generality, we set $\tau = 1$ and assume the entire state is the output, i.e. $C=I \in \RR^{n \times n}$.

Here we use the Exact Algorithm~\ref{alg:MCBiCG} with preconditioner $K^{-1}$ and the adjoint preconditioner $K^{-T}$. The application of the preconditioners is analogous to the implementation of $(K-\sigma M)^{-1}$ and $(K-\sigma M)^{-T}$ with $\sigma=0$. The shifted preconditioned system \eqref{shifted-c} is approximated with one execution of the Exact Algorithm~\ref{alg:MCBiCG}, generating approximations to $x(\mu)$ as in \eqref{our-problem}. The relative residual norm at iteration $i$ is computed as in \eqref{rel-res-pde} with $A(\mu)$~\eqref{time-delay-sys}. 

\begin{figure}
\centering
\includegraphics[]{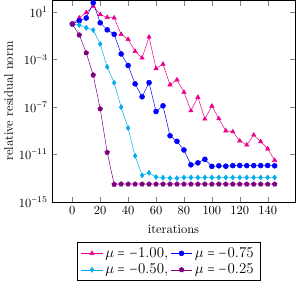}
\hspace{1.6cm}
\includegraphics[]{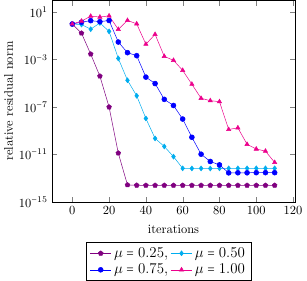}
\caption{Convergence for approximating \eqref{our-problem}, generated from one (total) execution of the Exact Algorithm~\ref{alg:MCBiCG} with $\sigma=0$ and $A(\mu)$ as in \eqref{time-delay-sys}. Here $d=17$, $tol = 10^{-11}$, relative residual norm \eqref{rel-res-pde}, and $b$ a random vector.}
\label{fig5}
\end{figure}

The results of this experiment are in Figure~\ref{fig5}. We see that our method is competitive for a variety of positive and negative values of $\mu$, and that approximations corresponding to values of $\mu$ closer to the target parameter $\sigma = 0$ converge faster than approximations corresponding to values of $\mu$ farther away. As the Exact Algorithm~\ref{alg:MCBiCG} is run, each approximation is updated using just additional scalar and vector computations. Before each approximate solution to \eqref{our-problem} is returned, a final application of the preconditioner is performed. The basis for the Krylov subspace \eqref{krylov-subspace} does not need to be stored if the values of the parameter $\mu$ are determined in advance. Additionally, the larger the dimension $j$ of the Krylov subspace from which the approximates come, the more solutions we have access to. As before, the nonlinear function $e^{-\mu}$ is approximated using \texttt{Chebfun} in Matlab with truncation parameter $d$ as in~\eqref{chebyshev}. 

\section{An inexact variant of preconditioned Chebyshev BiCG for parameterized linear systems} \label{sec:Flex}
In the short-term recurrence method the Exact Algorithm~\ref{alg:MCBiCG}, the action of the preconditioners $(K-\sigma M)^{-1}$ and $(K-\sigma M)^{-T}$ are applied via a block LU decomposition of the matrix product $(K-\sigma M) \Pi$ in \eqref{permutation-matrix} when approximating the solution to the shifted preconditioned linear system \eqref{shifted}. In this way, each application of $(K-\sigma M)^{-1}$ requires the action of $P(\sigma)^{-1} \in \mathbb{R}^{n \times n}$, and each application of $(K-\sigma M)^{-T}$ requires the action of $P(\sigma)^{-T} \in \mathbb{R}^{n \times n}$, as can be seen in \eqref{s-and-i} and \eqref{s-and-i-T}. An LU decomposition of the matrix $P(\sigma)$ was performed in the precomputing step and reused at each iteration, though other choices based on a factorization are possible. This approach is only suitable when an LU decomposition of a matrix of dimension $n \times n$ is feasible.

Consider now an inexact preconditioned BiCG method, where the action of $(K-\sigma M)^{-1}$ as well as the action of $(K-\sigma M)^{-T}$ are applied approximately at each iteration. Specifically, the corresponding linear systems with system matrices $P(\sigma)$ and $P(\sigma)^T$ are solved iteratively, and the accuracy of the solves varies from one iteration to the next. This approach, inspired by the work \cite{Vogel2007}, eliminates the need for an LU decomposition of $P(\sigma)$, offering an improvement in performance for approximating solutions to \eqref{our-problem}. A similar inner-outer approach was investigated in \cite{10.1007/3-540-28504-0_13}, where the effect of the error in inexact matrix-vector products was analyzed for several different iterative methods. Additionally, in \cite{KurschFreitag}, a relaxation strategy for low-rank ADI was studied for solving Lyapunov equations. Both of these works successfully increased the inner tolerance as the error in the outer method decreased, and we apply this strategy to our method. 

Similar approaches, where Krylov methods have been used as preconditioners, have been investigated in, for example, \cite{Chapman1998,deSturler1996,SimonciniSzyld2003}, as well as \cite{Axelsson1991,Szyld2001,Vuik19932}, where stopping criteria was utilized. To our knowledge, this is the first time an inexact BiCG method has been used to solve for the solution of multiple shifted systems simultaneously and the first time such a framework has been paired with a linearization of the form in \eqref{companion} to solve parameterized systems, where the dependence on the parameter is nonlinear.  

\subsection{Derivation of inexact preconditioned Chebyshev BiCG for parameterized linear systems} 
The $j$th iteration of the Lanczos biorthogonalization generated by the matrix $M (K-\sigma M)^{-1}$ and its adjoint appear in \eqref{relations}. As we are interested in an inexact algorithm, we consider the analogous relations corresponding to the inexact and iteration-dependent application of the preconditioning matrices $(K-\sigma M)^{-1}$ and $(K-\sigma M)^{-T}$, i.e., on iteration $j$ the relations
\begin{subequations} \label{flex-fac}
\begin{alignat}{5}
	M \hat{Z}_j &= V_j &&\hat{T}_j &&+ \hat{\beta}_j v_{j+1} e_j^T &&= V_{j+1} &&\underline{\hat{T}}_j ,\label{flex-fac-a} \\
	\hat{X}_j &= W_j &&\hat{T}_j^T &&+ \hat{\gamma}_j w_{j+1} e_j^T &&= W_{j+1} &&\hat{\bar{T}}^T_j
\end{alignat}
\end{subequations}
hold, where 
\begin{subequations} \label{comp-zs}
\begin{eqnarray}
\hat{z}_i &=& \mathcal{P}_{1_i}^{-1} v_i, \label{standard} \\
\hat{x}_i &=& \mathcal{P}_{2_i}^{-1} \left( M^T w_i \right)  \label{nonstandard}
\end{eqnarray}
\end{subequations}
with $\mathcal{P}_{1_i}^{-1} \approx (K-\sigma M)^{-1}$ and $\mathcal{P}_{2_i}^{-1} \approx (K-\sigma M)^{-T}$. 
Here $\hat{Z}_j = \begin{bmatrix} \hat{z}_1, \ldots, \hat{z}_j \end{bmatrix}$, $\hat{X}_j = \begin{bmatrix} \hat{x}_1, \ldots, \hat{x}_j \end{bmatrix} \in \mathbb{R}^{dn \times j}$, and the matrices $\hat{T}_j\in \mathbb{R}^{j \times j}$ and $\underline{\hat{T}}_j, \hat{\bar{T}}^T_j \in \mathbb{R}^{(j +1) \times j}$ are of the same form as the matrices in \eqref{square-T} and \eqref{Ts} with entries $\hat{\alpha}_i$, $\hat{\gamma}_i$, and $\hat{\beta}_i$, $i=1,\ldots,j$, defined below. Note that the choice of notation in \eqref{comp-zs} above is to differentiate the application of the preconditioner in \eqref{standard} from the application of the adjoint preconditioner in \eqref{nonstandard}. More specifically, the preconditioning matrix $(K-\sigma M)^{-T}$ is applied inexactly to the vector $(M^T w_i) \in \RR^{dn}$ in the adjoint case.  

In this setting, an application of the preconditioner $\mathcal{P}_{1_i}^{-1}$ refers to approximating the action of $P(\sigma)^{-1}$ within the application of $L_{\sigma}^{-1}$ as in \eqref{s-and-i}, and an application of $\mathcal{P}_{2_i}^{-1}$ refers to the application of $P(\sigma)^{-T}$ in the application of $L_{\sigma}^{-T}$ as in \eqref{s-and-i-T}. We denote an exact application of the preconditioners on the vectors $v_i$ and $M^T w_i$ as
\begin{subequations} \label{zs}
\begin{eqnarray}
	z_i &\coloneqq& (K-\sigma M)^{-1} v_i, \\
	x_i &\coloneqq& (K-\sigma M)^{-T} \left( M^T w_i \right).
\end{eqnarray}
\end{subequations}
Note that these vectors are not computed in the inexact algorithm, as the preconditioners are applied only approximately. These vectors are defined for the purpose of our analysis. 

Consider the $i$th column in equations \eqref{flex-fac}, i.e., 
\begin{subequations} \label{i-th-col}
\begin{alignat}{3} 
	M \hat{z}_i &= \hat{\gamma}_{i-1} v_{i-1} &&+ \hat{\alpha}_i v_i &&+ \hat{\beta}_i v_{i+1}, \\
	\hat{x}_i &= \hat{\beta}_{i-1} w_{i-1} &&+ \hat{\alpha}_i w_i &&+ \hat{\gamma}_i w_{i+1}, 
\end{alignat}
\end{subequations}
where $v_0 \coloneqq 0$, $w_0 \coloneqq 0$. Equations \eqref{i-th-col} paired with the biorthogonality condition \eqref{iden} lead to the definition
\begin{align} \label{alpha-i}
	\hat{\alpha}_i \coloneqq w_i^T M \hat{z}_i. 
\end{align}
Further, we define
\begin{subequations} \label{res}
\begin{alignat}{5} 
	\hat{r}_i &\coloneqq \hat{\beta}_i v_{i+1} &&= M &&\hat{z}_i &&- \hat{\alpha}_i v_i &&- \hat{\gamma}_{i-1} v_{i-1},\label{res-a} \\
	\hat{s}_i &\coloneqq \hat{\gamma}_i w_{i+1} &&= &&\hat{x}_i &&- \hat{\alpha}_i w_i &&- \hat{\beta}_{i-1} w_{i-1}, \label{res-b}
\end{alignat}
\end{subequations}
giving
\[
	1 = w_{i+1}^T v_{i+1} = \left( \frac{\hat{s}_i}{\hat{\gamma}_i} \right)^T \left( \frac{\hat{r}_i}{\hat{\beta}_i} \right),
\]
and thus the following definition: 
\begin{align} \label{gamma-i}
	\hat{\gamma}_i \coloneqq \frac{\hat{s}_i^T \hat{r}_i}{\hat{\beta}_i} \cdot
\end{align}
We define $\hat{\beta}_i$ using the so-called canonical choice as in \cite{GoluVanl96}, i.e., 
\begin{align} \label{beta-i}
	\hat{\beta}_i \coloneqq \norm{\hat{r}_i}_2.
\end{align}

We are interested in solving the shifted preconditioned system \eqref{shifted} with our inexact BiCG method. Equivalent to \eqref{flex-fac}, on iteration $j$ the shifted inexact relations
\begin{subequations} \label{krylov-rel}
\begin{alignat}{2}
	V_j + (-\mu + \sigma) M \hat{Z}_j &= V_{j+1} &&\Big(\underline{I}_j + (-\mu + \sigma) \underline{\hat{T}}_j \Big), \\
	W_j + (-\mu + \sigma) \hat{X}_j  &= W_{j+1} &&\Big(\underline{I}_j + (-\mu + \sigma) \hat{\bar{T}}^T_j \Big)
\end{alignat}
\end{subequations}
hold, where $\underline{I}_j$ is as in \eqref{shifted-relations}. In practice, we form the matrices $\hat{Z}_j$ and $\hat{\underline{T}}_j$ in \eqref{krylov-rel} once and compute $\tilde{x}_j (\mu)$, the approximate the solution to \eqref{our-problem} on iteration $j$, for each value of $\mu$ as 
\begin{subequations} \label{flex-approx}
\begin{eqnarray}
	y_j(\mu) &=& (I_j + (-\mu + \sigma) \hat{T}_j)^{-1} (\beta e_1), \label{flex-a} \\
	\hat{u}_j(\mu) &=& \hat{Z}_j y_j (\mu), \label{flex-b} \\
	\tilde{x}_j (\mu) &=& (\hat{u}_j(\mu))_{1:n},  \label{flex-c}
\end{eqnarray}
\end{subequations}
where $\beta \coloneqq \norm{b}$ and $(\hat{u}_j(\mu))_{1:n}$ denotes the first $n$ entries of $(\hat{u}_j(\mu)) \in \mathbb{R}^{dn}$. Here $I_j$ is as in~\eqref{iden}, $e_1$ is the first column of the identity matrix of dimension $j \times j$, and $\hat{T}_j \in \mathbb{R}^{j \times j}$ is defined as $\underline{\hat{T}}_j$ with the last row removed. We assume the matrix in \eqref{flex-a} is nonsingular.

Computing the approximation $\tilde{x}_j (\mu)$ requires one linear solve with a tridiagonal matrix of dimension $j \times j$ and one matrix-vector product with a matrix of dimension $dn \times j$ for each value of $\mu$. Note, if at any iteration the inner product of the residual vectors vanishes, i.e., $(\hat{s}_i)^T \hat{r}_i=0$, the algorithm has a breakdown. This situation never occurred while carrying out the experiments in this work. A full description of inexact preconditioned Chebyshev BiCG for parameterized linear systems appears in Algorithm~\ref{alg:FMCBiCG}. We will refer to Algorithm~\ref{alg:FMCBiCG} as the Inexact Algorithm~\ref{alg:FMCBiCG} or the {\em inexact algorithm}. 

The Inexact Algorithm~\ref{alg:FMCBiCG} terminates when the approximations $\tilde{x}_j(\mu_l)$ to the linear systems given in \eqref{our-problem} from the Krylov subspace of dimension $j$ have relative residual norm below a certain tolerance $tol$, for $l=1,\ldots,k$. As the Inexact Algorithm~\ref{alg:FMCBiCG} stores the basis matrix $\hat{Z}_j$ for approximations to $\tilde{x}_j (\mu)$, as well as the tridiagonal matrix $\hat{T}_j$, approximations to \eqref{our-problem} corresponding to $\mu \notin \{ \mu_l \},$ $l=1,\ldots,k$, can be computed after the algorithm has been executed once. In particular, it is reasonable to expect accurate approximations corresponding to values of $\mu$ close to the target $\sigma$, i.e., for $\mu$ such that $\abs{\mu - \sigma } \leq \abs{ \mu_k-\sigma }$. 

\begin{remark}[Storing $\hat{Z}_j$ and $\hat{T}_j$]
In the Exact Algorithm~\ref{alg:MCBiCG}, the colinearity of the residuals of the $i$th approximations of BiCG applied to the seed and shifted systems was used in order to derive a short-term recurrence method for shifted systems. These residual vectors spanned \eqref{krylov-subspace}, and a basis matrix for the Krylov subspace was not stored. Further, updates were computed based on an implicit LU factorization of $T_j$ in \eqref{relations}, i.e., by using a coupled two-term recurrence. This is also the approach in the standard BiCG method. 

It was not possible to show an analogous colinearity of residuals for the inexact method, so the Inexact Algorithm~\ref{alg:FMCBiCG} does not update approximations to the shifted systems with the same approach as the Exact Algorithm~\ref{alg:MCBiCG}. For reasons of presentation, the Inexact Algorithm~\ref{alg:FMCBiCG} stores the basis matrix $\hat{Z}_j$ and computes the approximations corresponding to each $\mu$ as described in \eqref{flex-approx}. An approach based on an implicit LU factorization of each matrix $(I_j + (-\mu_l + \sigma) \hat{T}_j)$, $l=1,\ldots,k$, could have been taken. Note, while $\hat{Z}_j$ is a basis matrix for the set of approximations given by the Inexact Algorithm~\ref{alg:FMCBiCG}, it is not a basis matrix for \eqref{krylov-subspace}; see Remark~\ref{not-krylov}. 
\end{remark}

\begin{remark}[Numerical behavior of the Exact Algorithm~\ref{alg:MCBiCG} and the Inexact Algorithm~\ref{alg:FMCBiCG}]
The Exact Algorithm~\ref{alg:MCBiCG} and the Inexact Algorithm~\ref{alg:FMCBiCG} with exact preconditioning are equivalent to applying the standard BiCG method to the shifted linear system \eqref{last-in-shifted} but are based on different approaches. Both methods return the same approximate solution, ignoring roundoff errors, but methods which use the same approach as in the Exact Algorithm~\ref{alg:MCBiCG}, i.e., ones based on a coupled two-term recurrence, are often more robust and have better numerical behavior; see \cite{FreundNachtigal}. 
\end{remark}

\begin{remark}[Krylov subspaces of the {\em inexact algorithm}] \label{not-krylov}
A result in \cite{inexKry} implies that the relations \eqref{flex-fac} can be expressed equivalently as
\begin{subequations} \label{flex-fac2}
\begin{alignat}{4}
	M \hat{Z}_j &=  \left(M\hphantom{t} (K-\sigma M)^{-1} \hphantom{t}+ \mathcal{E}_j \right) &&V_j &&= V_{j+1} && \underline{\hat{T}}_j , \\
	\hat{X}_j &=  \left( (K-\sigma M )^{-T} M^T +  \mathcal{F}_j \right) &&W_j &&= W_{j+1} && \hat{\bar{T}}^T_j,
\end{alignat}
\end{subequations}
where $\mathcal{E}_j = \sum_{i=1}^j E_i v_i w_i^T$ and $\mathcal{F}_j = \sum_{i=1}^j F_i w_i v_i^T$ with $E_i$, $F_i \in \mathbb{R}^{dn \times dn}$. Here $v_i$, $w_i$ are the $i$th columns of $V_j$ and $W_j$, respectively, and the biorthogonality condition \eqref{iden} holds. The matrices $E_i$ and $F_i$ represent the error which is introduced on each inexact application of the preconditioner and its adjoint. Note, the relations \eqref{flex-fac2} imply that the columns of $V_j$ in \eqref{flex-fac} span the Krylov subspace generated on the matrix $\left(M (K-\sigma M)^{-1} + \mathcal{E}_j \right)$ and the vector $v_1$, and, analogously, the columns of $W_j$ in \eqref{flex-fac} span the Krylov subspace generated from the matrix $\left( (K-\sigma M )^{-T} M^T + \mathcal{F}_j \right) $ and the vector $w_1$. We do not compute the matrices $E_i$ or $F_i$ and mention them purely for theoretical reasons. 
\end{remark}

\begin{algorithm}[h!]
\SetAlgoLined
\SetKwInput{Input}{Input}\SetKwInput{Output}{Output}\SetKwInput{Initialize}{Initialize}
\SetKwFor{FOR}{for}{do}{end}
\SetKwFor{IF}{if}{do}{end}
\SetKwFor{FUN}{function}{}{end}
\Input{$P_{\ell}$, $\ell = 0,1,\ldots,d$, as in \eqref{chebyshev} (Chebyshev coefficients) \\ $\tilde{b} \in \RR^{dn}$ as in \eqref{b-tilde}, $\tilde{c} \in \RR^{dn}$ such that $\tilde{b}^T \tilde{c} \neq 0$ \\ $\sigma \in \RR$ as in \eqref{shifted} (target parameter), $tol$ \\ $\{ \mu_l \}_{l=1,\ldots,k}$, $\mu_l \in \mathbb{R}$, $|\mu_1 - \sigma| \leq |\mu_2 - \sigma| \leq \cdots \leq |\mu_k - \sigma|$ (ordered set of shifts) \\ $A(\mu_l)$, $l=1,\ldots,k$, as in \eqref{summation}}
\Output{Approx. sol. $\tilde{x}_{j}(\mu_l)$, $l=1,\ldots,k$, to \eqref{our-problem}, from subspace of dim. $j$, matrices $\hat{Z}_j \in \mathbb{R}^{dn \times j}$, $\hat{T}_j\in\mathbb{R}^{j \times j}$}
\Initialize{$\hat{r}_0 = \tilde{b}$, $\hat{s}_0 = \tilde{c}$, $\hat{Z}_0 = []$}
\FOR{$i = 1,2, \dots$}{
Normalize $v_{i} = \hat{r}_{i-1} / \hat{\beta}_{i-1}$ using \eqref{beta-i}  \\
Normalize $w_i = \hat{s}_{i-1} / \hat{\gamma}_{i-1} $ using \eqref{gamma-i} \\
Compute $\hat{z}_i$ and $\hat{x}_i$ as in \eqref{comp-zs}\\
Update $\hat{Z}_i = \begin{bmatrix} \hat{Z}_{i-1} & \hat{z}_i \end{bmatrix}$ \\
Compute $\hat{\alpha}_i$ as in \eqref{alpha-i} \\
Compute $\hat{r}_i$ as in \eqref{res-a} \\
Compute $\hat{s}_i$ as in \eqref{res-b} \\
Update $\hat{T}_i$ as in \eqref{flex-fac} \\
$\tilde{x}_{i}(\mu_k)=$\textsc{ PostProcess}($\hat{T}_i,\hat{Z}_i,\mu_k,\sigma$) \\
$res = \norm{A(\mu_k) \tilde{x}_{i}(\mu_k) - b}/\norm{b}$ \\
\IF {$res \leq tol$}
{Set $\tilde{x}_{j}(\mu_k)=\tilde{x}_{i}(\mu_k)$, $j = i$ \\ \FOR{l=1,\ldots,k-1}{$\tilde{x}_{j}(\mu_l)=$\textsc{ PostProcess}($\hat{T}_j,\hat{Z}_j,\mu_l,\sigma$)} 
\IF 
{$\norm{A(\mu_l) \tilde{x}_{j}(\mu_l) - b}/\norm{b} \leq tol$, $l=1,\ldots,k-1$} 
{\textbf{return} } }
}
\FUN {$\tilde{x}_{j}(\mu)=$ \textsc{PostProcess}($\hat{T}_j,\hat{Z}_j,\mu,\sigma$)} {
Compute $y_{j}(\mu)$ as in \eqref{flex-a} \\
Compute $\hat{u}_j (\mu)$ as in \eqref{flex-b} \\
\textbf{return} $\tilde{x}_{j}(\mu)$ as in \eqref{flex-c}
}
\caption{Inexact preconditioned Chebyshev BiCG for parameterized linear systems}
\label{alg:FMCBiCG}
\end{algorithm}

\subsection{Explicit computation of the residual in the {\em inexact algorithm}}
The Inexact Algorithm~\ref{alg:FMCBiCG} approximates the solution to the linear system \eqref{companion} by solving the shifted preconditioned linear system \eqref{shifted}, where the preconditioners $(K- \sigma M)^{-1}$ and $(K - \sigma M)^{-T}$ are applied approximately in an inexact BiCG setting. In order to better understand the convergence of our approach, we compute a bound on the residual at each iteration of the {\em inexact algorithm} applied to this system. This bound includes a contribution that is directly related to the error in the application of the preconditioner. 

Define the inner residual vector $p_i$ on iteration $i$ as
\begin{align} \label{res-p} 
	p_i \coloneqq (K - \sigma M) \hat{z}_i - v_i,
\end{align}
where $\hat{z}_i$ is as in \eqref{comp-zs}. Equivalently, in matrix form, $\hat{Z}_j \coloneqq Z_j + (K-\sigma M)^{-1} P_j$, where $\hat{Z}_j$ is as in \eqref{flex-fac}, $Z_j = \begin{bmatrix} z_1, \ldots, z_j \end{bmatrix} \in \RR^{dn \times j}$ with $z_i$ as in \eqref{zs}, and $P_j  = \begin{bmatrix} p_1, \ldots, p_j \end{bmatrix} \in \RR^{dn \times j}$. The inner residual vector represents how inexactly we apply the preconditioner on the $i$th iteration of the Inexact Algorithm~\ref{alg:FMCBiCG}.

The relations \eqref{flex-fac} and \eqref{krylov-rel} are equivalent to the shifted relations given by
\begin{subequations} \label{krylov-rel2}
\begin{alignat}{4}
	V_j + (-\mu + \sigma) M \hat{Z}_j &= V_j &&\left(I_j + (-\mu + \sigma) \hat{T}_j \right) &&+ (-\mu+\sigma) \hat{\beta}_j v_{j+1} e_j^T, \\
	W_j + (-\mu + \sigma) \hat{X}_j  &= W_j &&\left(I_j + (-\mu + \sigma) \hat{T}_j^T\right) &&+ (-\mu+\sigma) \hat{\gamma}_j w_{j+1} e_j^T,
\end{alignat}
\end{subequations}
and $r_i^{in}$, the residual of the $i$th iterate of the Inexact Algorithm~\ref{alg:FMCBiCG} applied to the linear system~\eqref{companion}, is expressed as 
\begin{subequations} \label{res-fmcbicg}
\begin{eqnarray}
	r_i^{in} &=& \tilde{b} - (K-\mu M) \hat{u}_i(\mu) \\
	&=& \tilde{b} - (K-\mu M) \hat{Z}_i y_i (\mu) \\
	&=& \tilde{b} - (K-\mu M + \sigma M - \sigma M) \hat{Z}_i y_i(\mu) \\
	&=& \tilde{b} - ((-\mu + \sigma) M + (K - \sigma M)) \hat{Z}_i y_i(\mu) \\
	&=& \tilde{b} - (-\mu + \sigma) M \hat{Z}_i y_i(\mu) - V_i y_i(\mu) - P_i y_i(\mu) \\
	&=& \left( V_{i} \left( \beta e_1 - \left(I_i + (-\mu + \sigma) \hat{T}_i \right) \right) - (-\mu+\sigma) \hat{\beta}_i v_{i+1} e_i^T\right) y_i(\mu) - P_i y_i(\mu) \\
	&=& (\mu-\sigma) \hat{\beta}_i v_{i+1} e_i^T  y_i(\mu) - P_i y_i(\mu) \\
	&=& r^{app}_i - P_i y_i(\mu), \label{this-line}  
\end{eqnarray}
\end{subequations} 
where $\hat{u}_i(\mu)$, $y_i(\mu)$ are as in \eqref{flex-approx}, $\beta \coloneqq \norm{\tilde{b}}$ with $\tilde{b}$ is as in \eqref{linear}, and we have used the shifted relations \eqref{krylov-rel2}. The vector $r_i^{app}$ approximates $r_i^{ex}$, the residual of the Exact Algorithm~\ref{alg:MCBiCG} applied to \eqref{companion}. This residual is defined as $r_i^{ex} \coloneqq (\mu-\sigma) \beta_i v_{i+1} e_i^T  y_i(\mu)$, where $\beta_i$ and $v_{i+1}$, as well as the computation of $y_i(\mu)=(I_i + (-\mu + \sigma) T_i)^{-1} (\beta e_1)$, similar to \eqref{flex-approx}, stem from the relations \eqref{relations}; see \cite{IteratMSparseSys}. Note, the quantity $r_i^{app}$ above is efficient to compute. An analogous result regarding the application of the adjoint preconditioner $(K-\sigma M)^{-T}$ holds for the residual of the adjoint linear system given in \eqref{adjoint}.

\subsection{Convergence of the {\em inexact algorithm}} \label{sec:further-analysis}
Preconditioned Krylov subspace methods are only suitable when the action of the preconditioner is cheap to apply. It is, therefore, of interest to apply the preconditioners in the Inexact Algorithm~\ref{alg:FMCBiCG} in the most efficient way possible. Inspired by work in \cite{inexKry}, we prove a computable bound on the residual of the {\em inexact algorithm} applied to the system \eqref{shifted}. The error in the inexact application of the preconditioner contributes to this bound.

Consider first the shifted relations \eqref{krylov-rel} after $j-1$ iterations of the inexact Lanczos biorthogonalization process. Let a QR decomposition of the shifted tridiagonal upper Hessenberg matrix be such that 
\begin{subequations}\label{first-qr}
\begin{eqnarray} 
	\bm{\underline{T}}_{j-1} &\coloneqq&  \left( \underline{I}_{j-1} + (-\mu + \sigma) \underline{\hat{T}}_{j-1} \right) \\
	&=& Q_{j-1} \begin{bmatrix} R_{j-1}^T & 0 \end{bmatrix}^T \\ &=&  Q_{j-1} \hat{R}_{j-1},
\end{eqnarray}
\end{subequations}
where the entries of $\bm{\underline{T}}_{j-1} \in \RR^{j \times (j-1)}$ are denoted as $\bm{\underline{t}}_{l,\ell}$, $l=1,\ldots,j$, $\ell = 1,\ldots,j-1$, the matrix $Q_{j-1}^T$ is defined as
\begin{align} \label{Q-givens}
	Q_{j-1}^T \coloneqq \Omega_{j-1} \Omega_{j-2} \cdots \Omega_1 \in \mathbb{R}^{j \times j},
\end{align}
i.e., the product of Givens rotation matrices $\Omega_i$, $i=1,\ldots,j-1$, given by
\[
	\Omega_i \coloneqq 
	\begin{bmatrix}
	I_{i-1} & & & \\
	& c_i & s_i & \\
	& -s_i & c_i & \\
	& & & I_{j-1-i}
	\end{bmatrix} \in \mathbb{R}^{j \times j},
\]
and $R_{j-1} \in \RR^{(j-1) \times (j-1)}$ is upper triangular. Here $s_i$, $c_i$ are the sines and cosines of the Givens rotations constructed to eliminate the nonzero elements $\bm{\underline{t}}_{i+1,i}$ on the subdiagonal of $\bm{\underline{T}}_{j-1}$. 

On lines $10-14$ of the Inexact Algorithm~\ref{alg:FMCBiCG}, we update our approximation $\tilde{x}_j(\mu)$ as in~\eqref{flex-approx} by first performing a linear solve with the square matrix $\bm{{T}}_{j} \coloneqq (I_j + (-\mu + \sigma) \hat{T}_j) \in \mathbb{R}^{j \times j}$. Using the QR factorization in \eqref{first-qr}, the matrix $\bm{T}_j$ can be expressed as
\begin{subequations} \label{QR-of-Tj}
\begin{eqnarray} 
	\bm{T}_j &=& 
	\begin{bmatrix}
	\bm{\underline{T}}_{j-1} & \bm{t}_j
	\end{bmatrix} \\
	&=& 
	\begin{bmatrix}
	Q_{j-1} \hat{R}_{j-1} & \bm{t}_j
	\end{bmatrix} \\
	&=& Q_{j-1}
	\begin{bmatrix}
	\hat{R}_{j-1} & Q_{j-1}^T \bm{t}_j
	\end{bmatrix} \\
	&=&
	Q_{j-1} \tilde{R}_{j},
\end{eqnarray}
\end{subequations}
where $\bm{t}_j \in \mathbb{R}^j$ is the $j$th column of $\bm{T}_j$, $\tilde{R}_{j} \coloneqq \begin{bmatrix} \hat{R}_{j-1} & Q_{j-1}^T \bm{t}_j \end{bmatrix} \in \mathbb{R}^{j \times j}$, and $Q_{j-1}^T$ is as in~\eqref{Q-givens}. We rewrite the linear solve in \eqref{flex-a} as 
\begin{align} \label{compute-y-j-proof}
	y_j (\mu) = \bm{T}_j^{-1} (\beta e_1) = \tilde{R}_j^{-1} Q_{j-1}^T (\beta e_1)
\end{align}
and define $\hat{g}_j \coloneqq Q_{j-1}^T (\beta e_1) \in \RR^j$, giving $y_j (\mu) =\tilde{R}_j^{-1} \hat{g}_j$, where $\tilde{R}_j^{-1}$ is an upper triangular matrix as the inverse of an upper triangular matrix. As shown in \cite{IteratMSparseSys}, the entries of $\hat{g}_j$ are given by $\hat{g}_j = 
\begin{bmatrix}
\gamma_1 & \ldots & \gamma_j 
\end{bmatrix}^T$, where $\gamma_i \coloneqq |c_i s_1 \cdots s_{i-1} \beta |$, for $i=1,\ldots,j-1$, and the $j$th component is equal to
\begin{align} \label{gamma-j}
	\gamma_j \coloneqq |s_1 \cdots s_{j-1} \beta|.
\end{align}

Let $\eta_i^{(j)} = (y_j(\mu))_i$ denote the $i$th component of $y_j(\mu)$. Then, 
\begin{subequations} \label{bound-on-y}
\begin{eqnarray}
	|(y_j(\mu))_i| = | \eta_i^{(j)} |&=& |(\tilde{R}_j^{-1})_{\substack{i,1:j}}
	 \hat{g}_j| \\
	 &\leq& \norm{(\tilde{R}_j^{-1})_{\substack{i,i:j}}} \norm{ (\hat{g}_j)_{\substack{i:j}} } \\
	 &=& \norm{e_i^T \tilde{R}_j^{-1}} \norm{(\hat{g}_j)_{\substack{i:j}}} \\
	 &=& \norm{\tilde{R}_j^{-1}} \norm{(\hat{g}_j)_{\substack{i:j}}} \\ 
	 &=& \frac{1}{\sigma_j (\bm{T}_j)} \norm{(\hat{g}_j)_{\substack{i:j}}},
\end{eqnarray}
\end{subequations} 
 where $\sigma_j (\bm{T}_j)$ denotes the largest singular value of $\bm{T_j}$ and
 \begin{align*}
 	\norm{(\hat{g})_{\substack{i:j}}}^2 &= \gamma_i^2 + \gamma_{i+1}^2 + \cdots + \gamma_j^2 \\
	&= \beta^2 (|c_i s_1 s_2 \cdots s_{i-1} |^2 + |c_{i+1} s_1 s_2 \cdots s_i |^2 + \cdots + |s_1 s_2 \cdots s_{j-1}|^2 \\ 
	&= \beta^2 (s_1 \cdots s_{i-1})^2 \left( |c_i|^2 + |c_{i+1} s_i |^2 + \cdots + | s_i \cdots s_{j-1} |^2 \right).
 \end{align*}
The equality $\norm{\begin{bmatrix} c_i, c_{i+1} s_i, \ldots, s_i \cdots s_{j-1} \end{bmatrix}} = \norm{\Omega_{j-1} \cdots \Omega_{i} e_i} = 1$ holds, and thus,  
\begin{align} \label{g-comp}
	\norm{(\hat{g})_{\substack{i:j}}} = \beta | s_1 \cdots s_{i-1}|.
\end{align}
Note, the norm of $r_i^{app}$, as computed in \eqref{res-fmcbicg}, is given by
\begin{align} \label{r-k-ex-norm}
	\norm{ r_i^{app}} = \norm{(-\mu+\sigma) \hat{\beta}_i v_{i+1} e_i^T  y_i(\mu)} = \frac{|\bm{\underline{t}}_{i+1,i}|}{| \tilde{r}_{i,i}|} | s_1 \cdots s_{i-1} \beta|,
\end{align}
where $(-\mu + \sigma)\hat{\beta}_i = \bm{\underline{t}}_{i+1,i}$, $\tilde{r}_{i,i}$ denotes the entry in the $i$th row, $i$th column of $\tilde{R}_i \in \RR^{i \times i}$, and we have computed the $i$th entry of $y_i(\mu) \in \RR^i$ using \eqref{compute-y-j-proof} and \eqref{gamma-j}. This computation is analogous to \cite[equation (5.2)]{Brown91}. We define the quantity 
\begin{align} \label{delta-k}
	\Delta_i \coloneqq \frac{| \tilde{r}_{i,i} |}{|\bm{\underline{t}}_{i+1,i}|} \norm{ r_i^{app} } = |\tilde{r}_{i,i} ||(y_i(\mu))_i | = \beta | s_1 \cdots s_{i-1}|, 
\end{align}
and thus obtain the bound
\begin{align} \label{y-bound}
	|(y_j(\mu))_i|=|\eta_i^{(j)} | \leq \frac{1}{\sigma_j(\bm{T}_j)} \Delta_i
\end{align}
from \eqref{bound-on-y} and \eqref{delta-k}. The following theorem expresses a computable bound on the norm of the residual of the {\em inexact algorithm} applied to the system \eqref{shifted}.

\begin{theorem} \label{conv-theorem}
Let $r_i^{in}$ be the residual of the $i$th iterate of the Inexact Algorithm~\ref{alg:FMCBiCG} applied to the linear system \eqref{companion}, and define 
\[
	r_i^{app} \coloneq (\mu-\sigma) \hat{\beta}_i v_{i+1} e_i^T  y_i(\mu), 
\]
for $i=1,\ldots,j$, where $\hat{\beta}_i$, $v_{i+1}$ are as in \eqref{krylov-rel2} and $y_i(\mu) = \begin{bmatrix} \eta_1^{(i)}, \ldots, \eta_i^{(i)} \end{bmatrix} \in \mathbb{R}^i$ as in~\eqref{flex-approx}. If at each iteration $i \leq j$ the inner residual vector $p_i$ \eqref{res-p} satisfies 
\begin{align} \label{pi-bound}
	\norm{p_i} \leq \frac{1}{j} \frac{\sigma_j(\bm{T}_j)}{\Delta_i} \varepsilon \coloneqq \varepsilon_{\text{inner}}^{(i)},
\end{align}
with $\Delta_i$ \eqref{delta-k}, then 
\[
	\norm{r_i^{in}} \leq \norm{r_i^{app}} + \varepsilon. 
\]
\end{theorem}

{\em Proof.}
The proof follows directly from the reasoning above paired with \eqref{res-fmcbicg}, i.e., 
\begin{align*}
	\norm{r_j^{in}} &\leq \norm{r_i^{app}} + \norm{P_j y_j(\mu)} = \norm{r_i^{app}} + \norm{\sum_{i=1}^j p_i \eta_i^{(j)}} \leq \norm{r_i^{app}} + \sum_{i=1}^j \norm{p_i} |\eta_i^{(j)}| \leq \norm{r_i^{app}} + \varepsilon, 
\end{align*}
where we have used the bounds in \eqref{y-bound} and \eqref{pi-bound}. $\endproof$

\subsection{Approximation of a parameterized Helmholtz equation by the {\em inexact algorithm}} \label{sec:subsec-helmholtz2}
\begin{figure}[h!]
\centering
	\includegraphics[]{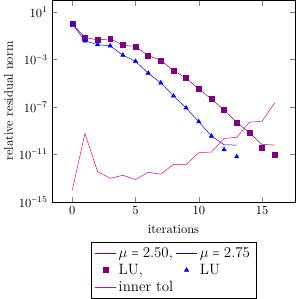}
	\hspace{1.6cm}
	\includegraphics[]{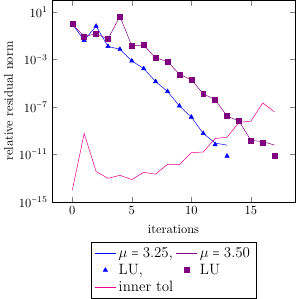}
\caption{Convergence for approximating the parameterized Helmholtz equation \eqref{pde}, generated from one (total) execution of the Inexact Algorithm~\ref{alg:FMCBiCG} with $\sigma=3$ and evaluated for different values of $\mu \in [-5,5]$. Here $n=243997$, $d=34$, (outer) $tol=10^{-10}$, and relative residual norm \eqref{rel-res-pde}. BiCG with variable inner tol for application of $P(\sigma)^{-1}$, according to \eqref{flex-tol-1} with $\varepsilon = 10^{-12}$. Compare to LU of $P(\sigma)$.}
\label{fig:flex}
\end{figure}

Consider approximating the parameterized Helmholtz equation~\eqref{pde} with the Inexact Algorithm~\ref{alg:FMCBiCG}. In the simulation shown in Figure~\ref{fig:flex}, the action of $P(\sigma)^{-1}$ in the application of the preconditioner is approximated via the iterative method BiCG. Here the tolerance in BiCG, referred to as the inner tolerance, is varied at each iteration. Specifically, we set the inner tolerance at iteration $i$, denoted $tol_i$, to $tol_1 = 10^{-14}$ and
\begin{align} \label{flex-tol-1}
	tol_i =  \frac{1} {|( (y_{i-1}(\mu^*))_{i-1}|} \varepsilon \approx \frac{1}{\Delta_{i-1}} \varepsilon, 
\end{align}
for $i=2,\ldots,j$, where $(y_{i-1}(\mu^*) )_{i-1}$ denotes the $(i-1)$th component of $y_{i-1}(\mu^*) \in \RR^{i-1}$ as in \eqref{flex-a} and $\Delta_{i-1}$ is as in \eqref{delta-k}. The parameter $\mu^*$ is equal to the $\mu_l$ furthest from the target parameter $\sigma$, i.e., the parameter with corresponding approximation from which we expect the slowest convergence (see \eqref{conv-reasoning}). Note, the inner tolerance as computed in \eqref{flex-tol-1} uses information from the previous iteration. 

\begin{figure}[h!]
\centering
	\includegraphics[]{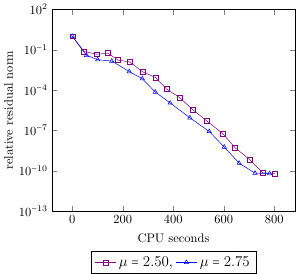}
	\hspace{1.6cm}
	\includegraphics[]{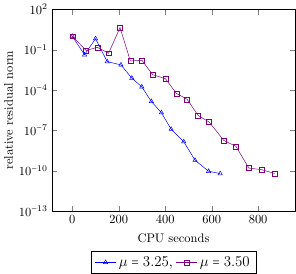}
\caption{Relative residual norm in terms of CPU sec for the simulation displayed in Figure~\ref{fig:flex}. Note, no LU decomposition of $P(\sigma)$.}
\label{fig:flex-timings}
\end{figure}

The experiment shown in Figure~\ref{fig:flex} was produced with one execution of the Inexact Algorithm~\ref{alg:FMCBiCG}. Though we display just four solutions, accessing each of the corresponding approximations to $\mu \in [2.50,3.50]$ requires the solution to a tridiagonal system of dimension $j \times j$, where $j$ is the dimension of the subspace from which the approximates come. The relative residual norms of these solutions are below a prescribed tolerance. Figure~\ref{fig:flex-timings} displays the convergence of the same simulation as a function of CPU time, omitting a comparison to an exact application of the preconditioner. In this way, we see the cost in CPU seconds of each iteration of the algorithm. The roughly constant cost of each iteration is due to the short-term recurrence feature of the method. The CPU times here are measured after the precomputation steps and, though one execution of the Inexact Algorithm~\ref{alg:FMCBiCG} gives approximations for many different values of $\mu$, each convergence curve here corresponds to a separate run. 

\begin{figure}[h!]
\centering
\begin{subfigure}[T]{0.46\columnwidth}
	\centering
	\includegraphics[]{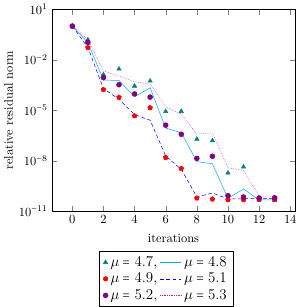}
	\caption{One execution of the Inexact Algorithm~\ref{alg:FMCBiCG} with $\sigma=5$ and (outer) $tol=10^{-10}$.}
	\label{fig:figlasta}
\end{subfigure}
\hspace{.5cm}
\begin{subfigure}[T]{0.46\columnwidth}
	\centering
	\includegraphics[]{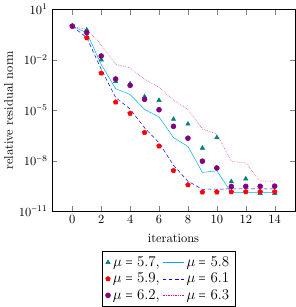}
\caption{One execution of the Inexact Algorithm~\ref{alg:FMCBiCG} with $\sigma=6$ and (outer) $tol=10^{-9}$.}
\label{fig:figlastb}
\end{subfigure}
\caption{Convergence for approximating the parameterized Helmholtz equation \eqref{pde}, generated from the Inexact Algorithm~\ref{alg:FMCBiCG} and evaluated for different values of $\mu \in [-8,8]$. Here $n=976076$, $d=44$, and relative residual norm \eqref{rel-res-pde}. AGMG with variable tol in application of $P(\sigma)^{-1}$ according to \eqref{flex-tol-1} with $\varepsilon = 10^{-12}$. Note, no LU of $P(\sigma)$ is feasible.}
\label{fig:flex2}
\end{figure}

Theorem~\ref{conv-theorem} shows a bound on the residual obtained from applying the {\em inexact algorithm} to the linearized system \eqref{shifted}, under the assumption that the bound in \eqref{pi-bound} is fulfilled for $i=1,\ldots,j$. As the quantity $\Delta_{i}$ is proportional to the $r_i^{app}$, an approximation to the residual assuming exact preconditioning, we can expect to solve the inner linear systems with growing inaccuracy as the outer residual decreases. Increasing the inner tolerance as the algorithm proceeds improves the performance without destroying the accuracy of the method. For comparison, we apply the preconditioner exactly via an LU decomposition of $P(\sigma)$ and display the convergence for a variety of values $\mu$. 

Figure~\ref{fig:flex2} shows the convergence of the {\em inexact algorithm} applied to the same discretization of the parameterized Helmholtz equation in \eqref{pde}, for larger values of the parameter $\mu$. Here we consider a problem of a larger dimension, where an LU decomposition of $P(\sigma)$ is no longer feasible, and we apply the action of $P(\sigma)^{-1}$ with the iterative method Aggregation-based algebraic multigrid\footnote{Yvan Notay, AGMG software and documentation; see \texttt{http://agmg.eu}} (AGMG) \cite{NapovNotay2012,Notay2010,Notay2012}. For this simulation, the inner tolerance is as described in \eqref{flex-tol-1}, i.e., we solve the inner systems with increasing inaccuracy as the outer residual of the method decreases. We see that for larger values of the parameter $\mu$, the target parameter $\sigma$ should be chosen closer to the values of $\mu$ which are of interest. As before, executing the algorithm once allows us to compute the corresponding approximates to all $\mu$ on a given interval in a cheap way, though we display just six of these solutions in the interval $[4.7,5.3]$ in Figure~\ref{fig:figlasta} and six solutions in $[5.7,6.3]$ in Figure~\ref{fig:figlastb}. 

In summary, the inexact framework successfully eliminates the need for an LU decomposition of the matrix $P(\sigma)$ in the precomputing steps, while maintaining the accuracy of the method for many values of $\mu$. The short-term recurrence property of the Inexact Algorithm~\ref{alg:FMCBiCG} offers a constant cost per iteration. Thus, the {\em inexact algorithm} is suitable for a wide range of large-scale simulations where an LU decomposition of a matrix is not feasible.

\section{Conclusions and future work} \label{sec:conclusion}
This work proposes two variants of a novel Krylov subspace method to approximate the solution to parameterized linear systems of the form~\eqref{our-problem}. Both algorithms return a function $\tilde{x}_j(\mu)$ on iteration $j$ which is cheap to evaluate for many different values of the parameter $\mu$. These algorithms are constructed by considering the approximate solution to a companion linearization based on an accurate Chebyshev interpolation of $A(\mu)$, where shift-and-invert preconditioners are used. The approximation to the resulting shifted preconditioned system is found in a shifted BiCG setting. 

Here both the preconditioner and its adjoint are applied via an efficient block LU decomposition of the matrix $(K-\sigma M) \Pi$ as in \eqref{permutation-matrix}. The first variant considers exact applications of the preconditioners, and the second variant applies an approximation to the preconditioners in an inexact setting. A computable bound on the residual obtained from iterates of the {\em inexact method} was shown, and a contribution in the bound is directly related to the error in the application of the preconditioner. Additionally, both algorithms offer a short-term recurrence, resulting in a constant cost per iteration. Numerical results confirm that the algorithms proposed here are suitable for large-scale simulations. 

The methods IDR(s) \cite{doi:10.1137/070685804} and IDR(s) for shifted systems \cite{DU201535} have proven effective for solving nonsymmetric linear systems. Another successful short-term recurrence method for shifted systems was developed in \cite{Frommer2003res}, based on the method Bi-CGSTAB \cite{doi:10.1137/0913035}. Using these methods to solve the linearization \eqref{shifted} would likely result in new robust methods, though further research would be needed. Furthermore, this work considered only real-valued preconditioners, though complex-valued shift-and-invert preconditioners have successfully been incorporated in several previous works for solving the Helmholtz equation; see, for instance, \cite{Erlangga08,ErlanggaEtAl2004}. While such a strategy would likely work here as well, it would require additional analysis.


\bibliographystyle{siam}
\bibliography{siobhanbib,eliasbib}

\end{document}